# Fast integration of DPG matrices based on tensorization

Jaime Mora[*] and Leszek Demkowicz

The Institute for Computational Engineering and Sciences (ICES), The University of Texas at Austin, 201 E 24th St, Austin, TX 78712, USA

## Abstract

Numerical integration of the stiffness matrix in higher order finite element (FE) methods is recognized as one of the heaviest computational tasks in a FE solver. The problem becomes even more relevant when computing the Gram matrix in the algorithm of the Discontinuous Petrov Galerkin (DPG) FE methodology. Making use of 3D tensor-product shape functions, and the concept of sum factorization, known from standard high order FE and spectral methods, here we take advantage of this idea for the entire exact sequence of FE spaces defined on the hexahedron. The key piece to the presented algorithms is the exact sequence for the one-dimensional element, and use of hierarchical shape functions. Consistent with existing results, the presented algorithms for the integration of $H^1$, $H(\text{curl})$, $H(\text{div})$, and $L^2$ inner products, have the $\mathcal{O}(p^7)$ computational complexity. Additionally, a modified version of the algorithms is proposed when the element map can be simplified, resulting in the reduced $\mathcal{O}(p^6)$ complexity. Use of Legendre polynomials for shape functions is critical in this implementation. Computational experiments performed with $H^1$, $H(\text{div})$ and $H(\text{curl})$ test shape functions show good correspondence with the expected rates.

## 1 Introduction

During the computation of a numerical solution to a linear boundary value problem using a Finite Element (FE) method, two major tasks are carried out by the computer processors, namely, calculation and assembly of the stiffness matrix and the load vector, and the solution of the linear system. The latter depends mostly on the size and properties of the global stiffness matrix, while the former depends more heavily on the local characteristics of every element in the mesh, as both the load vector and stiffness matrix are first computed at the element level. In higher-order finite elements the calculation of the integrals that compose such matrix and vector may become costly, since the higher the degree of a polynomial integrand, the larger the number of quadrature points needed to perform an exact (or at least a well approximated) numerical integration. Our interest in this article is to combine several results that can lead to significant savings in three-dimensional DPG computations.

In FEM, every type of element can be associated to one or more families of shape functions whether for the trial space or for the test space. In two spatial dimensions, the only conventional tensor-product element is the quadrilateral, which is generated by tensor-multiplying two line segments. In 3D, the conventional element types are the hexahedron, the tetrahedron, the (triangular) prism and the pyramid (of quadrilateral base) [12]. Out of them, the hexahedron and the prism are examples of tensor-product elements. The hexahedron is a triple tensor-product of 1D intervals, and the prism is a tensor-product of a triangle in 2D and a 1D interval. Thus, the standard shape functions for these two types of elements are tensor products of the shape functions associated to their lower-dimension generating elements.

---

[*]Corresponding author. E-mail: jmorapaz@ices.utexas.edu

The idea of tensorization for computing 2D and 3D integrals has been an established technique for spectral methods in CFD [22, 14], known as sum factorization, and later adopted for $p$- and $hp$-FEM with higher-order shape functions by Melenk et al [19]. It was additionally implemented for the case of a fully automatic $hp$-adaptive FEM solution of the Helmholtz equation by Kurtz [18], therein delivering explicit steps to compute the matrix while keeping a low memory requirement. More recently, the tensor-product nature of the shape functions in Bernstein-Bezier FEM and Isogeometric Analysis also motivated the extension of tensorized integration for the stiffness matrices of those methods [1, 2].

Use of fast integration algorithms becomes even more critical in an efficient implementation of the Discontinuous Petrov Galerkin (DPG) FE methodology. In order to see that, let us first refer to a classic Galerkin FE method in 3D. Let $p$ denote the order of polynomial basis for discretizing the solution. The algebraic structure of the final linear system reads

$$\mathsf{B}\mathsf{u} = \mathsf{l},$$

with square stiffness matrix $\mathsf{B}$ of size $\mathcal{O}(p^3) \times \mathcal{O}(p^3)$, solution vector $\mathsf{u}$ and load vector $\mathsf{l}$, both of the same size $\mathcal{O}(p^3)$. In DPG we enrich the test space usually by increasing the polynomial order in $\Delta p$ only for the test functions. Here, that single equation above is replaced by a larger system,

$$\begin{pmatrix} \mathsf{G} & \mathsf{B} & \widetilde{\mathsf{B}} \\ \mathsf{B}^\mathsf{T} & 0 & 0 \\ \widetilde{\mathsf{B}}^\mathsf{T} & 0 & 0 \end{pmatrix} \begin{pmatrix} \mathsf{s} \\ \mathsf{u} \\ \mathsf{w} \end{pmatrix} = \begin{pmatrix} \mathsf{l} \\ 0 \\ 0 \end{pmatrix}.$$

with additional unknowns $\mathsf{s}$ (size $\mathcal{O}((p+\Delta p)^3)$) and $\mathsf{w}$ (size $\mathcal{O}(p^2)$), and matrices $\mathsf{G}$ (square, size $\mathcal{O}(p+\Delta p)^3$), $\mathsf{B}$ (resized to $\mathcal{O}((p+\Delta p)^3) \times \mathcal{O}(p^3)$) and $\widetilde{\mathsf{B}}$ (size $\mathcal{O}((p+\Delta p)^3) \times \mathcal{O}(p^2)$), and load vector $\mathsf{l}$ (resized to $\mathcal{O}((p+\Delta p)^3)$). A conventional numerical integration of all these matrices has a complexity with leading terms of $\mathcal{O}((p+\Delta p)^9) + \mathcal{O}((p+\Delta p)^6 p^3)$ floating point operations. Furthermore, we can reduce this system by static condensation of $\mathsf{s}$, obtaining

$$\begin{pmatrix} \mathsf{B}^\mathsf{T}\mathsf{G}^{-1}\mathsf{B} & \mathsf{B}^\mathsf{T}\mathsf{G}^{-1}\widetilde{\mathsf{B}} \\ \widetilde{\mathsf{B}}^\mathsf{T}\mathsf{G}^{-1}\mathsf{B} & \widetilde{\mathsf{B}}^\mathsf{T}\mathsf{G}^{-1}\widetilde{\mathsf{B}} \end{pmatrix} \begin{pmatrix} \mathsf{u} \\ \mathsf{w} \end{pmatrix} = \begin{pmatrix} \mathsf{B}^\mathsf{T}\mathsf{G}^{-1}\mathsf{l} \\ \widetilde{\mathsf{B}}^\mathsf{T}\mathsf{G}^{-1}\mathsf{l} \end{pmatrix}.$$

Reaching this point involves a cost of $\mathcal{O}((p+\Delta p)^9)$ operations for the Cholesky factorization and necessary substitution steps to get $G^{-1}\mathsf{l}$, $\mathcal{O}((p+\Delta p)^9 p^3)$ for $G^{-1}\mathsf{B}$ and $G^{-1}\widetilde{\mathsf{B}}$, and a leading cost of $\mathcal{O}(p^6(p+\Delta p)^6)$ for the final matrix-matrix multiplications. These costs are almost unavoidable in order to have a solution using the DPG methodology, and we can handle them in an efficient way using highly specialized linear algebra libraries, which results in a smaller CPU time. Thus, if we want to get time savings in this technique, the part on which we can focus is the computation of $\mathsf{G}$, $\mathsf{B}$, $\widetilde{\mathsf{B}}$ and $\mathsf{l}$. Matrix $\mathsf{G}$, hereinafter referred to as the Gram matrix, is of special interest as it is the largest array involved. Therefore, having as a basis the sum factorization approaches mentioned earlier, we aim in this paper to adapt those algorithms to the computation of $\mathsf{G}$ and finally obtain a significant saving in the implementation of DPG. For more details on the derivation of the shown DPG system of equations see Section 3.

In the present paper, when deriving the algorithms we restrict ourselves to working with Gram matrices coming from the standard inner products only; however, it will be shown in one of the application cases that these ideas can also be extended to a more involved type of inner product. We also work just on the hexahedron case. For the subject of FE shape functions, this article uses as its main reference a thorough review of shape functions for elements of all shapes done by Fuentes et al [12]. We now outline the document: in Section 2 we explain the process of tensorization and develop algorithms to apply the technique for all the involved Hilbert spaces and simplify the operations under some assumptions. In Section 3 we provide several examples of DPG implementations and observe the computing time for the Gram matrix and other matrices, followed by a discussion. In Section 4 we close with conclusions and mention future extensions of the present work.

## 2 Tensorization

The approach sought in this article is related to the concept of exact sequences of energy spaces. With that starting point, we go from infinite-dimensional energy spaces to finite-dimensional subspaces which in the end are the ones we implement in any finite element method. We later define tensor-product spaces and proceed to study how to compute the Gram matrix for each space of the exact sequence for the hexahedron of the first type [21].

### 2.1 Exact sequences

Let $X_0$, $X_1$,..., $X_{N_s}$ be a finite family of vector spaces. Let $A_i : X_{i-1} \longrightarrow X_i$ be a linear operator, for $i = 1, ..., N_s$. We say that the following sequence or *complex*:

$$X_0 \xrightarrow{A_1} X_1 \xrightarrow{A_2} \cdots \xrightarrow{A_{N_s}} X_{N_s}$$

is exact if, for $i = 1, 2, ..., N_s - 1$, it holds that $\mathsf{R}(A_i) = \mathsf{N}(A_{i+1})$ (where $\mathsf{R}$ denotes the range of the operator and $\mathsf{N}$ symbolizes the nullspace or kernel). In the context of energy spaces, i.e. Hilbert spaces for the solution of variational formulations of boundary value problems, we will define for every number of spatial dimensions at least one exact sequence, specifying both the energy spaces and operators involved.

In the paper, $\Omega$ will denote a bounded and simply connected domain in $\mathbb{R}^N$, $N = 1, 2, 3$.

#### 2.1.1 Exact sequence in 1D ($\Omega = (a, b)$, $a, b \in \mathbb{R}$, $a < b$)

$$\mathbb{R} \xrightarrow{id} H^1(\Omega) \xrightarrow{\partial} L^2(\Omega) \xrightarrow{0} \{0\}$$

The presence of a zero operator in the last link of the exact sequence indicates that $\partial$ is a surjection. Similarly, the presence of the first link means that the nullspace of $\partial$ consists only of constant functions. Hereinafter we omit writing both the first and last links of the sequence, remembering that the last operator must be a surjection, and that the nullspace of the first one is made only by constant functions.

### 2.1.2 Exact sequences in 2D

In two dimensions, there are two possible exact sequences. They are:

$$H^1(\Omega) \xrightarrow{\nabla} H(\text{curl}, \Omega) \xrightarrow{\nabla_{vts} \times} L^2(\Omega),$$

$$H^1(\Omega) \xrightarrow{\nabla_{stv} \times} H(\text{div}, \Omega) \xrightarrow{\text{div}} L^2(\Omega),$$

where $\nabla_{vts} \times$ is the 2D vector-to-scalar curl operator and $\nabla_{stv} \times$ is the 2D scalar-to-vector curl operator. For a two-dimensional vector-valued function $E = (E_1, E_2) \in H(\text{curl}, \Omega)$, its curl is defined $\nabla_{vts} \times E = \partial_1 E_2 - \partial_2 E_1$; in turn, for a scalar $u \in H^1(\Omega)$ its curl is given by $\nabla_{stv} \times u = (\partial_2 u, -\partial_1 u)$.

Our work will focus on the three-dimensional problem therefore leaving the 2D case merely as a simplification of the discussed algorithms.

### 2.1.3 Exact sequence in 3D

Finally, the three-dimensional exact sequence involves all three classical vector calculus' differential operators,

$$H^1(\Omega) \xrightarrow{\nabla} H(\text{curl}, \Omega) \xrightarrow{\text{curl}} H(\text{div}, \Omega) \xrightarrow{\text{div}} L^2(\Omega). \tag{2.1}$$

We are interested in tensor-product three-dimensional finite elements. The most relevant example is the hexahedron, since it is a triple tensor product of the 1D simplex, the interval.

## 2.2 Tensor-product finite element shape functions

Let $\mathcal{I} = (0, 1)$ be the master interval in $\mathbb{R}$ and denote the master hexahedron by $\hat{\mathcal{K}} := \mathcal{I}^3$. Let $\mathbf{x}_\mathcal{K} : \hat{\mathcal{K}} \longrightarrow \mathcal{K}$ be the element map, transforming the master element into a physical space element $\mathcal{K}$. Suppose $\mathbf{x}_\mathcal{K}$ is a diffeomorphism over $\hat{\mathcal{K}}$. Then the Jacobian matrix $\mathcal{J}$ is well defined,

$$\mathcal{J} := \frac{\partial \mathbf{x}_\mathcal{K}}{\partial \boldsymbol{\xi}}, \tag{2.2}$$

where $\boldsymbol{\xi}$ is the position vector in the parametric (master) domain. The determinant of the Jacobian will be denoted $|\mathcal{J}| := \det \mathcal{J}$.

Let $T^{\text{grad}} : H^1(\hat{\mathcal{K}}) \longrightarrow H^1(\mathcal{K})$ be the map that takes the $H^1$ finite element space defined on the master domain $\hat{\mathcal{K}}$ to the physical space element $\mathcal{K}$. In the same fashion consider analogue maps denoted $T^{\text{curl}}$, $T^{\text{div}}$, $T$; these are the pullback or Piola maps for each function space in (2.1), and they are defined as follows:

$$H^1(\hat{\mathcal{K}}) \ni \hat{u} \longmapsto T^{\text{grad}} \hat{u} := \hat{u} \circ \mathbf{x}_\mathcal{K}^{-1} = u \in H^1(\mathcal{K}) \tag{2.3}$$

$$H(\text{curl}, \hat{\mathcal{K}}) \ni \hat{E} \longmapsto T^{\text{curl}} \hat{E} := (\mathcal{J}^{-T} \hat{E}) \circ \mathbf{x}_\mathcal{K}^{-1} = E \in H(\text{curl}, \mathcal{K}) \tag{2.4}$$

$$H(\text{div}, \hat{\mathcal{K}}) \ni \hat{V} \longmapsto T^{\text{div}} \hat{V} := (|\mathcal{J}|^{-1} \mathcal{J} \hat{V}) \circ \mathbf{x}_\mathcal{K}^{-1} = V \in H(\text{div}, \mathcal{K}) \tag{2.5}$$

$$L^2(\hat{\mathcal{K}}) \ni \hat{q} \longmapsto T \hat{q} := (|\mathcal{J}|^{-1} \hat{q}) \circ \mathbf{x}_\mathcal{K}^{-1} = q \in L^2(\mathcal{K}) \tag{2.6}$$

As an important remark concerning notation, we have opted to use a circumflex or "hat" ( ˆ ) above nearly every function, domain or vector space defined in the master space.

Now, it is known that for a FE method implementation, we don't work with the entire energy space but with a finite-dimensional linear subspace, usually consisting of polynomials. Following the notation in [10, 12], let us call such subspaces as follows:

$$\hat{W}^p \subsetneq H^1(\hat{\mathcal{K}})$$
$$\hat{Q}^p \subsetneq H(\operatorname{curl}, \hat{\mathcal{K}})$$
$$\hat{V}^p \subsetneq H(\operatorname{div}, \hat{\mathcal{K}})$$
$$\hat{Y}^p \subsetneq L^2(\hat{\mathcal{K}})$$

where the superindex $p$ symbolizes the nominal polynomial order of the sequence of spaces. Making sure that these finite-dimensional subspaces form themselves an exact sequence is important as it leads to useful properties when studying interpolants and approximability. We take that information as given because proving that fact is not central to this paper, thereby we just refer to those proofs within [10]. Through the Piola transformations we can get the resulting finite-dimensional subspaces in the physical space element. Those would be:

$$\begin{aligned} T^{\operatorname{grad}} \hat{W}^p =: W^p &\subsetneq H^1(\mathcal{K}) \\ T^{\operatorname{curl}} \hat{Q}^p =: Q^p &\subsetneq H(\operatorname{curl}, \mathcal{K}) \\ T^{\operatorname{div}} \hat{V}^p =: V^p &\subsetneq H(\operatorname{div}, \mathcal{K}) \\ T \hat{Y}^p =: Y^p &\subsetneq L^2(\mathcal{K}) \end{aligned} \quad (2.7)$$

Due to the element map being a diffeomorphism, there is a unique correspondence between any function of the physical element subspaces and the master element subspaces:

$$\forall u \in W^p \ \exists! \ \hat{u} \in \hat{W}^p \ \text{such that} \ u = T^{\operatorname{grad}} \hat{u} \tag{2.8}$$
$$\forall E \in Q^p \ \exists! \ \hat{E} \in \hat{Q}^p \ \text{such that} \ E = T^{\operatorname{curl}} \hat{E} \tag{2.9}$$
$$\forall V \in V^p \ \exists! \ \hat{V} \in \hat{V}^p \ \text{such that} \ V = T^{\operatorname{div}} \hat{V} \tag{2.10}$$
$$\forall q \in Y^p \ \exists! \ \hat{q} \in \hat{Y}^p \ \text{such that} \ q = T \hat{q}. \tag{2.11}$$

Those finite-dimensional subspaces for the hexahedron are [21]:

$$\hat{W}^p = \mathcal{Q}^{p_1, p_2, p_3}(\mathcal{I}^3)$$
$$\downarrow \nabla$$
$$\hat{Q}^p = \mathcal{Q}^{p_1-1, p_2, p_3}(\mathcal{I}^3) \times \mathcal{Q}^{p_1, p_2-1, p_3}(\mathcal{I}^3) \times \mathcal{Q}^{p_1, p_2, p_3-1}(\mathcal{I}^3)$$
$$\downarrow \operatorname{curl}$$
$$\hat{V}^p = \mathcal{Q}^{p_1, p_2-1, p_3-1}(\mathcal{I}^3) \times \mathcal{Q}^{p_1-1, p_2, p_3-1}(\mathcal{I}^3) \times \mathcal{Q}^{p_1-1, p_2-1, p_3}(\mathcal{I}^3)$$
$$\downarrow \operatorname{div}$$
$$\hat{Y}^p = \mathcal{Q}^{p_1-1, p_2-1, p_3-1}(\mathcal{I}^3)$$

where $p_1, p_2, p_3$ are positive integers and for any non-negative integers $p, q, r$ the space $\mathcal{Q}^{p,q,r}(\mathcal{I}^3) := \mathcal{P}^p(\mathcal{I}) \otimes \mathcal{P}^q(\mathcal{I}) \otimes \mathcal{P}^r(\mathcal{I})$, with $\mathcal{P}^p(\mathcal{I})$ being the space of univariate polynomials defined over $\mathcal{I}$ of degree less than or equal to $p$. A space like $\mathcal{Q}^{p,q,r}(\mathcal{I}^3)$ is known as a tensor product polynomial space, and an element of such a type of space is next characterized:

$$f \in \mathcal{Q}^{p,q,r}(\mathcal{I}^3) \iff \exists! f_1 \in \mathcal{P}^p(\mathcal{I}), f_2 \in \mathcal{P}^q(\mathcal{I}), f_3 \in \mathcal{P}^r(\mathcal{I}) \text{ such that } f(\boldsymbol{\xi}) = f_1(\xi_1)f_2(\xi_2)f_3(\xi_3) \; \forall \boldsymbol{\xi} \in \mathcal{I}^3. \tag{2.12}$$

Corresponding polynomial subspaces are also defined for the one-dimensional interval's exact sequence. They are:

$$\hat{W}_\mathcal{I}^p = \mathcal{P}^p(\mathcal{I}) \subsetneq H^1(\mathcal{I})$$
$$\downarrow \partial$$
$$\hat{Y}_\mathcal{I}^p = \mathcal{P}^{p-1}(\mathcal{I}) \subsetneq L^2(\mathcal{I})$$

Following the property of the exact sequence explained above, it must hold that $\partial \hat{W}_\mathcal{I}^p = \hat{Y}_\mathcal{I}^p$, which is easy to verify.

Notice how the hexahedron's exact sequence may be reconstructed by using the 1D interval's exact sequence multiple times.

$$\hat{W}^p = \hat{W}_\mathcal{I}^{p_1} \otimes \hat{W}_\mathcal{I}^{p_2} \otimes \hat{W}_\mathcal{I}^{p_3}$$
$$\downarrow \nabla = (\partial_1, \partial_2, \partial_3)$$
$$\hat{Q}^p = \hat{Y}_\mathcal{I}^{p_1} \otimes \hat{W}_\mathcal{I}^{p_2} \otimes \hat{W}_\mathcal{I}^{p_3} \times \hat{W}_\mathcal{I}^{p_1} \otimes \hat{Y}_\mathcal{I}^{p_2} \otimes \hat{W}_\mathcal{I}^{p_3} \times \hat{W}_\mathcal{I}^{p_1} \otimes \hat{W}_\mathcal{I}^{p_2} \otimes \hat{Y}_\mathcal{I}^{p_3}$$
$$\downarrow \mathrm{curl} = (\partial_2(\cdot)_3 - \partial_3(\cdot)_2, \partial_3(\cdot)_1 - \partial_1(\cdot)_3, \partial_1(\cdot)_2 - \partial_2(\cdot)_1)$$
$$\hat{V}^p = \hat{W}_\mathcal{I}^{p_1} \otimes \hat{Y}_\mathcal{I}^{p_2} \otimes \hat{Y}_\mathcal{I}^{p_3} \times \hat{Y}_\mathcal{I}^{p_1} \otimes \hat{W}_\mathcal{I}^{p_2} \otimes \hat{Y}_\mathcal{I}^{p_3} \times \hat{Y}_\mathcal{I}^{p_1} \otimes \hat{Y}_\mathcal{I}^{p_2} \otimes \hat{W}_\mathcal{I}^{p_3}$$
$$\downarrow \mathrm{div} = \partial_1(\cdot)_1 + \partial_2(\cdot)_2 + \partial_3(\cdot)_3$$
$$\hat{Y}^p = \hat{Y}_\mathcal{I}^{p_1} \otimes \hat{Y}_\mathcal{I}^{p_2} \otimes \hat{Y}_\mathcal{I}^{p_3}$$

This shows that, if in a FE element subroutine we replace the calls to 3D shape functions by multiple calls to the 1D shape functions, we can reconstruct all the spaces in the exact sequence. As it will be explained below, doing this can provide great benefits with regard to computational performance

Finally, the Gram matrix that is going to be studied throughout this work is explicitly defined as follows. Let $\mathcal{H}$ be a finite-dimensional Hilbert space, with inner product $(\cdot, \cdot)_\mathcal{H}$, $\dim \mathcal{H} = N_h$, and a basis $\{h_i\}_{i=1}^{N_h}$. The Gram matrix $G^\mathcal{H}$ is given by:

$$G_{ij}^\mathcal{H} := (h_i, h_j)_\mathcal{H} \text{ for } i, j = 1, ..., N_h. \tag{2.13}$$

Our goal is to propose algorithms to compute a Gram matrix (2.13) more efficiently than the conventional ones when dealing with tensor-product finite-element spaces of shape functions. All the energy spaces in the 3D exact sequence (2.1) are going to be analyzed (setting $\Omega = \mathcal{K}$), thereby a particular inner product is required to be defined in every case. The finite-dimensional Hilbert spaces will be those defined in (2.7).

The following subsections are presented in order of complexity instead of their position in the exact sequence. Unless otherwise specified, the master element domain $\hat{\mathcal{K}}$ corresponds to $\mathcal{I}^3$. We also restrict ourselves to real-valued functions for simplicity, but no major modification is anticipated when moving into complex-valued functions.

## 2.3 Space $L^2$

Let us recall the definition of the $L^2(\mathcal{K})$ space:

$$L^2(\mathcal{K}) = \left\{ \text{Lebesgue-measurable functions } f : \mathcal{K} \to \mathbb{R} \ : \int_{\mathcal{K}} |f(\boldsymbol{x})|^2 d^3\boldsymbol{x} < \infty \right\}.$$

The symbol for the inner product in $L^2(\mathcal{K})$, $(\cdot, \cdot)_{L^2(\mathcal{K})}$, typically incorporates the domain of integration as a subscript:

$$(\varphi, \vartheta)_{\mathcal{K}} := (\varphi, \vartheta)_{L^2(\mathcal{K})} = \int_{\mathcal{K}} \varphi(\boldsymbol{x})\vartheta(\boldsymbol{x})d^3\boldsymbol{x} \quad \forall\, \varphi, \vartheta \in L^2(\mathcal{K}). \tag{2.14}$$

For vector-valued functions living in $\boldsymbol{L}^2(\mathcal{K}) := (L^2(\mathcal{K}))^3$ the associated inner-product is defined componentwise, that is,

$$(\Phi, \Theta)_{\mathcal{K}} := (\Phi, \Theta)_{\boldsymbol{L}^2(\mathcal{K})} = \int_{\mathcal{K}} \Phi(\boldsymbol{x})^\mathsf{T} \Theta(\boldsymbol{x}) d^3\boldsymbol{x} = \sum_{d=1}^{3} (\varphi_d, \vartheta_d), \ \forall\, \Phi = (\varphi_1, \varphi_2, \varphi_3), \Theta = (\vartheta_1, \vartheta_2, \vartheta_3) \in \boldsymbol{L}^2(\mathcal{K}). \tag{2.15}$$

Let the order of the shape functions for the master hexahedron be $(p_1, p_2, p_3)$, in the sense of the exact sequence. Consider a basis for $Y^p$, $\{v_I\}_{I=0}^{\dim Y^p - 1}$, where $\dim Y^p = p_1 p_2 p_3$. Thus, for any pair of integers $0 \leq I, J < \dim Y^p$ the corresponding entry for the $L^2$ Gram matrix $G$ is obtained as follows:

$$\begin{aligned}
G_{IJ} &= (v_I, v_J)_{\mathcal{K}} \\
&= \int_{\mathcal{K}} v_I(\boldsymbol{x}) v_J(\boldsymbol{x}) d^3\boldsymbol{x} \\
&= \int_{\hat{\mathcal{K}}} v_I \circ \mathbf{x}_{\mathcal{K}}(\boldsymbol{\xi}) v_J \circ \mathbf{x}_{\mathcal{K}}(\boldsymbol{\xi}) |\mathcal{J}(\boldsymbol{\xi})| d^3\boldsymbol{\xi} \\
&= \int_{\hat{\mathcal{K}}} (T\hat{v}_I) \circ \mathbf{x}_{\mathcal{K}}(\boldsymbol{\xi}) (T\hat{v}_J) \circ \mathbf{x}_{\mathcal{K}}(\boldsymbol{\xi}) |\mathcal{J}(\boldsymbol{\xi})| d^3\boldsymbol{\xi} \\
&= \int_{\hat{\mathcal{K}}} [(|\mathcal{J}|^{-1} \hat{v}_I) \circ \mathbf{x}_{\mathcal{K}}^{-1}] \circ \mathbf{x}_{\mathcal{K}}(\boldsymbol{\xi}) [(|\mathcal{J}|^{-1} \hat{v}_J) \circ \mathbf{x}_{\mathcal{K}}^{-1}] \circ \mathbf{x}_{\mathcal{K}}(\boldsymbol{\xi}) |\mathcal{J}(\boldsymbol{\xi})| d^3\boldsymbol{\xi} \\
&= \int_{\hat{\mathcal{K}}} \hat{v}_I(\boldsymbol{\xi}) \hat{v}_J(\boldsymbol{\xi}) |\mathcal{J}(\boldsymbol{\xi})|^{-1} d^3\boldsymbol{\xi} \\
&= \int_0^1 \int_0^1 \int_0^1 \hat{v}_I(\xi_1, \xi_2, \xi_3) \hat{v}_J(\xi_1, \xi_2, \xi_3) |\mathcal{J}(\xi_1, \xi_2, \xi_3)|^{-1} d\xi_3 d\xi_2 d\xi_1,
\end{aligned} \tag{2.16}$$

where the definition of the $L^2$ Piola map (see (2.6)) was applied, and the element map between the master and physical space was invoked to transform the integrand. In the last line of the derivation (2.16) the volume integral over the master hexahedron is rewritten as three univariate integrals over the master interval. Now, since $\hat{v}_I$ and $\hat{v}_J$ belong to $\hat{Y}^p$, they are tensor-product polynomials in $\mathcal{Q}^{p_1-1,p_2-1,p_3-1}(\mathcal{I}^3)$; thus, if we take $\hat{v}_I$ as the model case, we have:

$$\hat{v}_I(\xi_1, \xi_2, \xi_3) := \nu_{1;i_1}(\xi_1) \nu_{2;i_2}(\xi_2) \nu_{3;i_3}(\xi_3) \tag{2.17}$$

where the univariate polynomials $\{\nu_{a;i_a}\}_{i_a=0}^{p_a-1}$ form a basis of shape functions for the space $\mathcal{P}^{p_a-1}(\mathcal{I})$, for $a = 1, 2, 3$, and the integer indices $0 \leq i_a < p_a$ are given so that they uniquely correspond to the original index $I$ (e.g. through the formula $I = i_1 + p_1 i_2 + p_1 p_2 i_3$). It is important to remark that should we account for a hierarchical basis of polynomials, then we could use that basis for each $\mathcal{P}^{p_a-1}(\mathcal{I})$, and the need for the first identifier in the index of $\nu_{a;i_a}$ goes away. Assuming that's the case, we can rewrite $\hat{v}_I$ and $\hat{v}_J$ as follows:

$$\hat{v}_I(\xi_1, \xi_2, \xi_3) = \nu_{i_1}(\xi_1)\nu_{i_2}(\xi_2)\nu_{i_3}(\xi_3)$$
$$\hat{v}_J(\xi_1, \xi_2, \xi_3) = \nu_{j_1}(\xi_1)\nu_{j_2}(\xi_2)\nu_{j_3}(\xi_3). \quad (2.18)$$

Combining (2.16) and (2.18), we get the following:

$$(v_I, v_J)_{\mathcal{K}} = \int_0^1 \int_0^1 \int_0^1 \nu_{i_1}(\xi_1)\nu_{i_2}(\xi_2)\nu_{i_3}(\xi_3)\nu_{j_1}(\xi_1)\nu_{j_2}(\xi_2)\nu_{j_3}(\xi_3)|\mathcal{J}(\xi_1, \xi_2, \xi_3)|^{-1} d\xi_3 d\xi_2 d\xi_1$$
$$= \int_0^1 \nu_{i_1}(\xi_1)\nu_{j_1}(\xi_1) \left\{ \int_0^1 \nu_{i_2}(\xi_2)\nu_{j_2}(\xi_2) \left[ \int_0^1 \nu_{i_3}(\xi_3)\nu_{j_3}(\xi_3)|\mathcal{J}(\xi_1, \xi_2, \xi_3)|^{-1} d\xi_3 \right] d\xi_2 \right\} d\xi_1 \quad (2.19)$$

The second line above depicts how the original volume integral turns into three nested interval integrals, through Fubini's theorem. Let us define a sequence of auxiliary functions with which we will assemble the triple integral in (2.19).

$$\mathcal{G}^A_{i_3 j_3}(\xi_1, \xi_2) := \int_0^1 \nu_{i_3}(\xi_3)\nu_{j_3}(\xi_3)|\mathcal{J}(\xi_1, \xi_2, \xi_3)|^{-1} d\xi_3 \quad (2.20)$$

$$\mathcal{G}^B_{i_2 j_2 i_3 j_3}(\xi_1) := \int_0^1 \nu_{i_2}(\xi_2)\nu_{j_2}(\xi_2)\mathcal{G}^A_{i_3 j_3}(\xi_1, \xi_2) d\xi_2 \quad (2.21)$$

$$\Rightarrow G_{IJ} = \mathcal{G}_{i_1 j_1 i_2 j_2 i_3 j_3} := \int_0^1 \nu_{i_1}(\xi_1)\nu_{j_1}(\xi_1)\mathcal{G}^B_{i_2 j_2 i_3 j_3}(\xi_1) d\xi_1 \quad (2.22)$$

Given $\xi_1, \xi_2$ we can evaluate $\mathcal{G}^A_{i_3 j_3}$ with numerical integration. Let $\xi_3^n$, $n = 1, ..., N$ be the collection of quadrature points in the interval (0,1) that make a polynomial integrand of degree $2N - 1$ be numerically integrated with full accuracy, and $w_3^N$ be the associated weight. Then the value of (2.20) may be approximated by:

$$\mathcal{G}^A_{i_3 j_3}(\xi_1, \xi_2) \approx \sum_{n=1}^N \nu_{i_3}(\xi_3^n)\nu_{j_3}(\xi_3^n)|\mathcal{J}(\xi_1, \xi_2, \xi_3^n)|^{-1} w_3^n. \quad (2.23)$$

In the same manner, let $\xi_2^m, w_2^m$, with $m = 1, ..., M$ and $\xi_1^l, w_1^l$, with $l = 1, ..., L$. We can approximate thus expressions (2.21) and (2.22) as follows:

$$\mathcal{G}^B_{i_1 j_2 i_3 j_3}(\xi_1) \approx \sum_{m=1}^M \nu_{i_2}(\xi_2^m)\nu_{j_2}(\xi_2^m)\mathcal{G}^A_{i_3 j_3}(\xi_1, \xi_2^m) w_2^m, \quad (2.24)$$

$$\mathcal{G}_{i_1 j_1 i_2 j_2 i_3 j_3} \approx \sum_{l=1}^L \nu_{i_1}(\xi_1^l)\nu_{j_1}(\xi_1^l)\mathcal{G}^B_{i_2 j_2 i_3 j_3}(\xi_1^l) w_1^l. \quad (2.25)$$

Making distinction of different sets of quadrature points makes sense when the individual polynomial degrees are not equal, but even in that case, things may get simpler if we choose to work only with the largest of those sets in all cases and we can therefore establish $L = M = N$, and $\xi_1^l = \xi_2^l = \xi_3^l = \zeta^l$, and the weight $w_1^l = w_2^l = w_3^l = w^l$ for every $l$ between 1 and $L$. Below it will be clear that even after taking this shortcut, the complexity of the resulting algorithm will not be driven by such a quadrature order.

Suppose we want to approximate the triple integral in the last line of (2.16) without taking advantage of the tensorization process. Then, we will have $L^3$ quadrature points, which are the triplets $(\xi_1^l, \xi_2^m, \xi_3^n) =: \boldsymbol{\xi}_{lmn}$, with their associated weights $w_{lmn} := w^l w^m w^n$ so the approximate inner product will be:

$$G_{IJ} = (v_I, v_J)_{\mathcal{K}} \approx \sum_{l=1}^{L} \sum_{m=1}^{L} \sum_{n=1}^{L} \hat{v}_I(\boldsymbol{\xi}_{lmn}) \hat{v}_J(\boldsymbol{\xi}_{lmn}) |\mathcal{J}(\boldsymbol{\xi}_{lmn})|^{-1} w_{lmn}. \tag{2.26}$$

Clearly, (2.26) and (2.23)-(2.25) are equivalent if we use the same quadrature order per coordinate. However, here lies the entire spirit of the sum-factorization or tensorization-based integration, as we will see. Firstly, it is a direct observation that the algorithm to compute (2.26) is the conventional one in 3D finite element codes, which is presented in Algorithm 1. Here, as usual in FE algorithms, the symmetry $G_{IJ} = G_{JI}$ is taken advantage of, so that the off-diagonal entries are computed only once.

---

**Algorithm 1** Conventional computation of the $L^2$ Gram Matrix

    **procedure** L2GRAM($mid, G$)     ▷ Compute matrix $G$ for element No. $mid$ - Conventional algorithm
        **call** setquadrature3D($mid, p_1, p_2, p_3, L, \boldsymbol{\xi}_{lmn}, w_{lmn}$)
        $G \leftarrow 0$     ▷ Initialize Gram Matrix
        **for** $l, m, n = 1$ to $L$ **do**
            **call** Shape3L2($\boldsymbol{\xi}_{lmn}, p_1, p_2, p_3, \{\hat{v}_I(\boldsymbol{\xi}_{lmn})\}$)     ▷ Evaluate 3D shape functions at $\boldsymbol{\xi}_{lmn}$
            **call** geometry( $\boldsymbol{\xi}_{lmn}, mid, \mathbf{x}, \mathcal{J}(\boldsymbol{\xi}_{lmn}), \mathcal{J}^{-1}(\boldsymbol{\xi}_{lmn}), |\mathcal{J}|$)     ▷ Compute jacobian
            **for** $J = 0$ to $\dim Y^p - 1$ **do**
                **for** $I = J$ to $\dim Y^p - 1$ **do**
                    $G_{IJ} \leftarrow G_{IJ} + \hat{v}_I(\boldsymbol{\xi}_{lmn}) \hat{v}_J(\boldsymbol{\xi}_{lmn}) |\mathcal{J}|^{-1} w_{lmn}$     ▷ Accumulate through (2.26)
    **return** $G$

---

If we have a uniform polynomial degree $p = p_1 = p_2 = p_3$, we will need at least $L = p$ to compute the approximate integral. Consequently, in Algorithm 1 the accumulation statement is going to be executed $\frac{1}{2}p^6(p^3 + 1)$ times. This represents an operation count of $\mathcal{O}(p^9)$.

Now, let us study the algorithm for the tensorization-based integration. The way we nested the interval integrals in (2.19) suggests we can perform a similar ordering in the algorithm loops. The idea is to fix a quadrature point in the coordinate 1, then compute the corresponding term in the sum of (2.25), for which we need to go over each quadrature point in coordinate 2 in order to obtain the sum in (2.24), and at each step of those it is required to go over all the quadrature points in the third coordinate and evaluate the expression (2.23); all must iterate until the sum to approximate the inner product is completed.

Additionally, see (2.20)-(2.22) to notice the symmetry in $\mathcal{G}^A_{i_3 j_3} = \mathcal{G}^A_{j_3 i_3}$ and the two minor symmetries in the other auxiliary function $\mathcal{G}^B_{i_2 j_2 i_3 j_3} = \mathcal{G}^B_{j_2 i_2 i_3 j_3} = \mathcal{G}^B_{i_2 j_2 j_3 i_3} = \mathcal{G}^B_{j_2 i_2 j_3 i_3}$. Furthermore we have the symmetry of the inner product, also noticeable in (2.22). Those symmetries can be exploited in order to make a more efficient computation of the Gram Matrix. The Algorithm 2 includes all the considerations above. Notice at the very end, the presence of a subroutine named $fillsymmetric3L2$, which has the goal of filling in (to the same extent than in the conventional algorithm) the Gram matrix $G$ by making use of the symmetries of the auxiliary functions $\mathcal{G}, \mathcal{G}^A, \mathcal{G}^B$. Of course, this function can be disregarded if the rest of the FE code is designed to work with $\mathcal{G}$ as an output.

---

**Algorithm 2** Tensorization-based computation of the $L^2$ Gram Matrix

---

  **procedure** L2GRAMTENSOR($mid, G$)    ▷ Compute matrix $G$ for element No. $mid$ - Tensorization-based algorithm
    $p_{\max} \leftarrow \max\{p_1, p_2, p_3\}$
    **call** setquadrature1D($mid, p_{\max} - 1, L, \zeta^l, w^l$)
    $\mathcal{G} \leftarrow 0$      ▷ Initialize Gram Matrix
    **for** $l = 1$ to $L$ **do**
      **call** Shape1L2($\zeta^l, p_1, \nu_{i_1}(\zeta^l)\}$)    ▷ Evaluate 1D shape functions at $\zeta^l$
      **for** $j_3 = 0$ to $p_3 - 1$ **do**
        **for** $i_3 = j_3$ to $p_3 - 1$ **do**
          $\mathcal{G}^B \leftarrow 0$
          **for** $m = 1, L$ **do**
            **call** Shape1L2($\zeta^m, p_2, \nu_{i_2}(\zeta^m)\}$)    ▷ Evaluate 1D shape functions at $\zeta^m$
            $\mathcal{G}^A \leftarrow 0$
            **for** $n = 1$ to $L$ **do**
              **call** Shape1L2($\zeta^n, p_3, \nu_{i_3}(\zeta^n)\}$)    ▷ Evaluate 1D shape functions at $\zeta^n$
              $\boldsymbol{\xi}_{lmn} \leftarrow (\zeta^l, \zeta^m, \zeta^n)$
              **call** geometry($\boldsymbol{\xi}_{lmn}, mid, \mathbf{x}, \mathcal{J}(\boldsymbol{\xi}_{lmn}), \mathcal{J}^{-1}(\boldsymbol{\xi}_{lmn}), |\mathcal{J}|$)    ▷ Compute jacobian
              $\mathcal{G}^A_{i_3 j_3} \leftarrow \mathcal{G}^A_{i_3 j_3} + \nu_{i_3}(\zeta^n)\nu_{j_3}(\zeta^n)|\mathcal{J}|^{-1}w^n$    ▷ Accumulate through (2.23)
            **for** $j_2 = 0$ to $p_2 - 1$ **do**
              **for** $i_2 = j_2$ to $p_2 - 1$ **do**
                $\mathcal{G}^B_{i_2 j_2 i_3 j_3} \leftarrow \mathcal{G}^B_{i_2 j_2 i_3 j_3} + \nu_{i_2}(\zeta^m)\nu_{j_2}(\zeta^m)\mathcal{G}^A_{i_3 j_3}(\zeta^l, \zeta^m)w^m$    ▷ From (2.24)
          **for** $j_2 = 0$ to $p_2 - 1$ **do**
            **for** $i_2 = j_2$ to $p_2 - 1$ **do**
              **for** $j_1 = 0$ to $p_1 - 1$ **do**
                **for** $i_1 = j_1$ to $p_1 - 1$ **do**
                  $\mathcal{G}_{i_1 j_1 i_2 j_2 i_3 j_3} \leftarrow \mathcal{G}_{i_1 j_1 i_2 j_2 i_3 j_3} + \nu_{i_1}(\zeta^l)\nu_{j_1}(\zeta^l)\mathcal{G}^B_{i_2 j_2 i_3 j_3}(\zeta^l)w^l$    ▷ From (2.25)
    **call** fillsymmetric3L2($p_1, p_2, p_3, \mathcal{G}, G$)    ▷ Fills matrix using symmetries
    **return** $G$

---

Assuming that Algorithm 2 is used having uniform polynomial degree $p$, then the cost of the integration is driven by the accumulation for each index combination of $\mathcal{G}_{i_1 i j_1 i_2 j_2 i_3 j_3}$, which is executed $p[{}^1\!/\!_2 p(p+1)]^3$ times, making a cost of $\mathcal{O}(p^7)$.

It is worth mentioning that the array $\mathcal{G}_{i_1 i j_1 i_2 j_2 i_3 j_3}$ needs not be constructed, especially if taking into account computer memory issues. Instead, its value can be directly accumulated into the corresponding entry of $G$.

In case the polynomial degrees are not uniform, different 1D rules could be applied per coordinate. If this is the case, a key point in the implementation of Algorithm 2 would be to reorder the integral nesting to make the quadrature points in the outermost loop correspond to the coordinate with the least polynomial order.

This saving of two orders of magnitude will be observed in all the remaining cases, although some may involve more complicated calculations.

**Remark 2.1.** The position of the loops for $i3, j3$ in Algorithm 2 was adopted from the one in Jason Kurtz's dissertation [18], and allows to save the memory associated to the indices $i3, j3$ in all the auxiliary arrays. We insist in this ordering for all the algorithms presented on this document because in DPG, as in any higher order FE technique, the polynomial degrees used could make the size of those arrays considerably large.

**Remark 2.2.** Algorithm 2 leads to an operation count of $\mathcal{O}(p^5)$ for the statement "call geometry(...)". Notice that for low values of $p$, this expensive call may cost a portion similar to the final accumulation statement. Bringing into consideration Remark 2.1, if memory limitations are not a major concern, an alternative way of implementing the sum factorization algorithm is shown in Algorithm 3, which reduces by two orders of magnitude the number of times "call geometry(...)" is executed.

## 2.4 Space $H^1$

This energy space is defined as follows:

$$H^1(\mathcal{K}) = \left\{ u \in L^2(\mathcal{K}) \; : \nabla u \in \boldsymbol{L}^2(\mathcal{K}) \right\} \tag{2.27}$$

The inner product in $H^1(\mathcal{K})$ can be computed by means of the expression:

$$(u, v)_{H^1(\mathcal{K})} := (u, v)_\mathcal{K} + (\nabla u, \nabla v)_\mathcal{K} \quad \forall u, v \in H^1(\mathcal{K}) \tag{2.28}$$

The gradient $\nabla$ used above is computed in the physical space, that is:

$$\nabla u = \left( \frac{\partial u}{\partial x_1}, \frac{\partial u}{\partial x_2}, \frac{\partial u}{\partial x_3} \right) = (\partial_1 u, \partial_2 u, \partial_3 u) \tag{2.29}$$

Moreover, we will need the gradient in the parametric space coordinates of a function $\hat{v}$ defined over $\hat{\mathcal{K}}$:

$$\hat{\nabla}\hat{v} = \left( \frac{\partial \hat{v}}{\partial \xi_1}, \frac{\partial \hat{v}}{\partial \xi_2}, \frac{\partial \hat{v}}{\partial \xi_3} \right) = \left( \hat{\partial}_1 \hat{v}, \hat{\partial}_2 \hat{v}, \hat{\partial}_3 \hat{v} \right) \tag{2.30}$$

The exact sequence in 3D described earlier, along with the Piola map definitions, imply that if $u \in H^1(\mathcal{K})$, then $\nabla u \in H(\text{curl}, \mathcal{K})$ and that if $u = T^{\text{grad}}\hat{u}$ for some $\hat{u} \in H^h at \mathcal{K})$, then $\nabla u = T^{\text{curl}}\hat{\nabla}\hat{u}$. Of course, also $\hat{\nabla}\hat{u} \in H(\text{curl}, \hat{\mathcal{K}})$.

Now, let the order of the shape functions for the master hexahedron be $(p_1, p_2, p_3)$, in the sense of the exact sequence. Consider a basis for $W^p$, $\{\varphi_I\}_{I=0}^{\dim W^p - 1}$, where $\dim W^p = (p_1+1)(p_2+2)(p_3+1)$. Thus, for any pair of integers $0 \leq I, J < \dim W^p$ the inner product is obtained as derived in (2.31).

Analogously to the $L^2$ case, besides recalling the definitions of both $T^{\text{grad}}$ and $T^{\text{curl}}$, in (2.31) we take the same steps as above to achieve an expression suitable for tensorization for computing the $H^1$ Gram matrix, herein denoted $G^{\text{grad}}$.

**Algorithm 3** Alternative tensorization-based computation of the $L^2$ Gram Matrix

**procedure** L2GRAMTENSOR($mid, G$) ▷ Compute matrix $G$ for element No. $mid$ - Tensorization-based algorithm
    $p_{\max} \leftarrow \max\{p_1, p_2, p_3\}$
    **call** setquadrature1D($mid, p_{\max} - 1, L, \zeta^l, w^l$)
    $\mathcal{G} \leftarrow 0$ ▷ Initialize Gram Matrix
    **for** $l = 1$ to $L$ **do**
        **call** Shape1L2($\zeta^l, p_1, \nu_{i_1}(\zeta^l)\}$) ▷ Evaluate 1D shape functions at $\zeta^l$
        $\mathcal{G}^B \leftarrow 0$
        **for** $m = 1, L$ **do**
            **call** Shape1L2($\zeta^m, p_2, \nu_{i_2}(\zeta^m)\}$) ▷ Evaluate 1D shape functions at $\zeta^m$
            $\mathcal{G}^A \leftarrow 0$
            **for** $n = 1$ to $L$ **do**
                **call** Shape1L2($\zeta^n, p_3, \nu_{i_3}(\zeta^n)\}$) ▷ Evaluate 1D shape functions at $\zeta^n$
                **for** $j_3 = 0$ to $p_3 - 1$ **do**
                      **for** $i_3 = j_3$ to $p_3 - 1$ **do**
                          $\boldsymbol{\xi}_{lmn} \leftarrow (\zeta^l, \zeta^m, \zeta^n)$
                          **call** geometry($\boldsymbol{\xi}_{lmn}, mid, \mathbf{x}, \mathcal{J}(\boldsymbol{\xi}_{lmn}), \mathcal{J}^{-1}(\boldsymbol{\xi}_{lmn}), |\mathcal{J}|$) ▷ Compute jacobian
                          $\mathcal{G}^A_{i_3 j_3} \leftarrow \mathcal{G}^A_{i_3 j_3} + \nu_{i_3}(\zeta^n)\nu_{j_3}(\zeta^n)|\mathcal{J}|^{-1} w^n$ ▷ Accumulate through (2.23)
            **for** $j_3 = 0$ to $p_3 - 1$ **do**
                **for** $i_3 = j_3$ to $p_3 - 1$ **do**
                      **for** $j_2 = 0$ to $p_2 - 1$ **do**
                          **for** $i_2 = j_2$ to $p_2 - 1$ **do**
                              $\mathcal{G}^B_{i_2 j_2 i_3 j_3} \leftarrow \mathcal{G}^B_{i_2 j_2 i_3 j_3} + \nu_{i_2}(\zeta^m)\nu_{j_2}(\zeta^m)\mathcal{G}^A_{i_3 j_3}(\zeta^l, \zeta^m) w^m$ ▷ From (2.24)
    **for** $j_3 = 0$ to $p_3 - 1$ **do**
        **for** $i_3 = j_3$ to $p_3 - 1$ **do**
            **for** $j_2 = 0$ to $p_2 - 1$ **do**
                **for** $i_2 = j_2$ to $p_2 - 1$ **do**
                      **for** $j_1 = 0$ to $p_1 - 1$ **do**
                          **for** $i_1 = j_1$ to $p_1 - 1$ **do**
                              $\mathcal{G}_{i_1 j_1 i_2 j_2 i_3 j_3} \leftarrow \mathcal{G}_{i j_1 i_2 j_2 i_3 j_3} + \nu_{i_1}(\zeta^l)\nu_{j_1}(\zeta^l)\mathcal{G}^B_{i_2 j_2 i_3 j_3}(\zeta^l) w^l$ ▷ From (2.25)
    **call** fillsymmetric3L2($p_1, p_2, p_3, \mathcal{G}, G$) ▷ Fills matrix using symmetries
    **return** $G$

$$\begin{aligned}
G^{\text{grad}}_{IJ} &= (\varphi_I, \varphi_J)_{H^1(\mathcal{K})} \\
&= \int_{\mathcal{K}} \varphi_I(\boldsymbol{x})\varphi_J(\boldsymbol{x}) d^3\boldsymbol{x} + \int_{\mathcal{K}} [\nabla\varphi_I(\boldsymbol{x})]^{\mathsf{T}} \nabla\varphi_J(\boldsymbol{x}) d^3\boldsymbol{x} \\
&= \int_{\hat{\mathcal{K}}} \varphi_I \circ \mathbf{x}_{\mathcal{K}}(\boldsymbol{\xi})\varphi_J \circ \mathbf{x}_{\mathcal{K}}(\boldsymbol{\xi})|\mathcal{J}(\boldsymbol{\xi})| d^3\boldsymbol{\xi} + \int_{\hat{\mathcal{K}}} [\nabla\varphi_I \circ \mathbf{x}_{\mathcal{K}}(\boldsymbol{\xi})]^{\mathsf{T}} \nabla\varphi_J \circ \mathbf{x}_{\mathcal{K}}(\boldsymbol{\xi})|\mathcal{J}(\boldsymbol{\xi})| d^3\boldsymbol{\xi}
\end{aligned}$$

$$\begin{aligned}
&= \int_{\hat{\mathcal{K}}} \left(T^{\text{grad}}\hat{\varphi}_I \circ \mathbf{x}_\mathcal{K}(\boldsymbol{\xi})\right)\left(T^{\text{grad}}\hat{\varphi}_J \circ \mathbf{x}_\mathcal{K}(\boldsymbol{\xi})\right)|\mathcal{J}(\boldsymbol{\xi})|d^3\boldsymbol{\xi} \; + \\
&\quad \int_{\mathcal{K}} \left[\nabla\left(T^{\text{grad}}\hat{\varphi}_I\right) \circ \mathbf{x}_\mathcal{K}(\boldsymbol{\xi})\right]^\mathsf{T}\left[\nabla\left(T^{\text{grad}}\hat{\varphi}_J\right) \circ \mathbf{x}_\mathcal{K}(\boldsymbol{\xi})\right]|\mathcal{J}(\boldsymbol{\xi})|d^3\boldsymbol{\xi} \\
&= \int_{\hat{\mathcal{K}}} \left((\hat{\varphi}_I \circ \mathbf{x}_\mathcal{K}^{-1}) \circ \mathbf{x}_\mathcal{K}(\boldsymbol{\xi})\right)\left((\hat{\varphi}_J \circ \mathbf{x}_\mathcal{K}^{-1}) \circ \mathbf{x}_\mathcal{K}(\boldsymbol{\xi})\right)|\mathcal{J}(\boldsymbol{\xi})|d^3\boldsymbol{\xi} \; + \\
&\quad \int_{\mathcal{K}} \left[\left(T^{\text{curl}}\hat{\nabla}\hat{\varphi}_I\right) \circ \mathbf{x}_\mathcal{K}(\boldsymbol{\xi})\right]^\mathsf{T}\left[\left(T^{\text{curl}}\hat{\nabla}\hat{\varphi}_J\right) \circ \mathbf{x}_\mathcal{K}(\boldsymbol{\xi})\right]|\mathcal{J}(\boldsymbol{\xi})|d^3\boldsymbol{\xi} \\
&= \int_{\hat{\mathcal{K}}} \hat{\varphi}_I(\boldsymbol{\xi})\hat{\varphi}_J(\boldsymbol{\xi})|\mathcal{J}(\boldsymbol{\xi})|d^3\boldsymbol{\xi} \; + \\
&\quad \int_{\mathcal{K}} \left[\left((\mathcal{J}^{-T}\hat{\nabla}\hat{\varphi}_I) \circ \mathbf{x}_\mathcal{K}^{-1}\right) \circ \mathbf{x}_\mathcal{K}(\boldsymbol{\xi})\right]^\mathsf{T}\left[\left((\mathcal{J}^{-T}\hat{\nabla}\hat{\varphi}_J) \circ \mathbf{x}_\mathcal{K}^{-1}\right) \circ \mathbf{x}_\mathcal{K}(\boldsymbol{\xi})\right]|\mathcal{J}(\boldsymbol{\xi})|d^3\boldsymbol{\xi} \\
&= \int_{\hat{\mathcal{K}}} \hat{\varphi}_I(\boldsymbol{\xi})\hat{\varphi}_J(\boldsymbol{\xi})|\mathcal{J}(\boldsymbol{\xi})|d^3\boldsymbol{\xi} + \int_\mathcal{K} \left[\hat{\nabla}\hat{\varphi}_I(\boldsymbol{\xi})\right]^\mathsf{T}\mathcal{D}(\boldsymbol{\xi})\left[\hat{\nabla}\hat{\varphi}_J(\boldsymbol{\xi})\right]|\mathcal{J}(\boldsymbol{\xi})|d^3\boldsymbol{\xi} \\
&= \int_0^1\int_0^1\int_0^1 \left\{\hat{\varphi}_I(\boldsymbol{\xi})\hat{\varphi}_J(\boldsymbol{\xi}) + \left[\hat{\nabla}\hat{\varphi}_I(\boldsymbol{\xi})\right]^\mathsf{T}\mathcal{D}(\boldsymbol{\xi})\left[\hat{\nabla}\hat{\varphi}_J(\boldsymbol{\xi})\right]\right\}|\mathcal{J}(\boldsymbol{\xi})|d\xi_3 d\xi_2 d\xi_1 \quad (2.31)
\end{aligned}$$

where the symmetric matrix $\mathcal{D} := \mathcal{J}^{-1}\mathcal{J}^{-T}$ contains the following entries:

$$\mathcal{D} = \begin{pmatrix} D_{11} & D_{12} & D_{13} \\ D_{21} & D_{22} & D_{23} \\ D_{31} & D_{32} & D_{33} \end{pmatrix}, \quad (2.32)$$

with $D_{12} = D_{21}$, $D_{13} = D_{31}$, $D_{32} = D_{23}$, and all of the entries being dependent on $\boldsymbol{\xi}$.

The conventional algorithm for (2.31) is presented as Algorithm 4.

---

**Algorithm 4** Conventional computation of the $H^1$ Gram Matrix

---

**procedure** H1GRAM($mid, G^{\text{grad}}$) ▷ Compute matrix $G^{\text{grad}}$ for element No. $mid$ - Conventional algorithm
    **call** setquadrature3D($mid, p_1, p_2, p_3, L, \boldsymbol{\xi}_{lmn}, w_{lmn}$)
    $G^{\text{grad}} \leftarrow 0$         ▷ Initialize Gram Matrix
    **for** $l, m, n = 1$ to $L$ **do**
        **call** Shape3H1($\boldsymbol{\xi}_{lmn}, p_1, p_2, p_3, \{\hat{\varphi}_I(\boldsymbol{\xi}_{lmn})\}, \{\hat{\nabla}\hat{\varphi}_I(\boldsymbol{\xi}_{lmn})\}$)   ▷ 3D shape functions at $\boldsymbol{\xi}_{lmn}$
        **call** geometry($\boldsymbol{\xi}_{lmn}, mid, \mathbf{x}, \mathcal{J}(\boldsymbol{\xi}_{lmn}), \mathcal{J}^{-1}(\boldsymbol{\xi}_{lmn}), |\mathcal{J}|$)   ▷ Compute jacobian
        $\mathcal{D} \leftarrow \mathcal{J}^{-1}(\boldsymbol{\xi}_{lmn})\mathcal{J}^{-T}(\boldsymbol{\xi}_{lmn})$
        **for** $J = 0$ to $\dim W^p - 1$ **do**
            **for** $I = J$ to $\dim W^p - 1$ **do**
                $G^{\text{grad}}_{IJ} \leftarrow G^{\text{grad}}_{IJ} + \left\{\hat{\varphi}_I(\boldsymbol{\xi}_{lmn})\hat{\varphi}_J(\boldsymbol{\xi}_{lmn}) + \left[\hat{\nabla}\hat{\varphi}_I(\boldsymbol{\xi}_{lmn})\right]^\mathsf{T}\mathcal{D}\left[\hat{\nabla}\hat{\varphi}_J(\boldsymbol{\xi}_{lmn})\right]\right\}|\mathcal{J}|w_{lmn}$
    **return** $G^{\text{grad}}$

---

The tensorization-based integration for the present situation begins with defining a tensor-product shape function. As $\hat{\varphi}_I, \hat{\varphi}_J \in \hat{W}^p = \mathcal{Q}^{p_1,p_2,p_3}(\mathcal{I}^3)$, it follows that:

$$\begin{aligned} \hat{\varphi}_I(\xi_1,\xi_2,\xi_3) &= \chi_{i_1}(\xi_1)\chi_{i_2}(\xi_2)\chi_{i_3}(\xi_3) \\ \hat{\varphi}_J(\xi_1,\xi_2,\xi_3) &= \chi_{j_1}(\xi_1)\chi_{j_2}(\xi_2)\chi_{j_3}(\xi_3). \end{aligned} \quad (2.33)$$

where the univariate polynomials $\{\chi_{i_a}\}_{i_a=0}^{p_a}$ are a hierarchical basis of the polynomial space $\mathcal{P}^{p_a}(\mathcal{I})$, for $a = 1, 2, 3$, and the integer indices $0 \leq i_a, j_a \leq p_a$ are given so that they hold a unique correspondence to the original indices $I, J$ (such as $I = i_1 + (p_1+1)i_2 + (p_1+1)(p_2+1)i_3$). Notice that the application of the master-domain gradient to the shape functions makes effect to a single 1D basis function per component, that is:

$$\hat{\nabla}\hat{\varphi}_I = \begin{pmatrix} \chi'_{i_1}(\xi_1)\chi_{i_2}(\xi_2)\chi_{i_3}(\xi_3) \\ \chi_{i_1}(\xi_1)\chi'_{i_2}(\xi_2)\chi_{i_3}(\xi_3) \\ \chi_{i_1}(\xi_1)\chi_{i_2}(\xi_2)\chi'_{i_3}(\xi_3) \end{pmatrix} \quad (2.34)$$

In the calculation of $G_{IJ}^{\text{grad}}$ the integrand term $\left[\hat{\nabla}\hat{\varphi}_I(\boldsymbol{\xi}_{lmn})\right]^\mathsf{T} \mathcal{D}\left[\hat{\nabla}\hat{\varphi}_J(\boldsymbol{\xi}_{lmn})\right]$ represents the main change with respect to the $L^2$ problem. Using (2.34) for the expanded form of that expression we derive the following (to see intermediate steps see [10]):

$$\begin{aligned} \left[\hat{\nabla}\hat{\varphi}_I\right]^\mathsf{T} \mathcal{D}\left[\hat{\nabla}\hat{\varphi}_J\right] &= \begin{pmatrix} \chi'_{i_1}\chi_{i_2}\chi_{i_3} & \chi_{i_1}\chi'_{i_2}\chi_{i_3} & \chi_{i_1}\chi_{i_2}\chi'_{i_3} \end{pmatrix} \begin{pmatrix} D_{11} & D_{12} & D_{13} \\ D_{21} & D_{22} & D_{23} \\ D_{31} & D_{32} & D_{33} \end{pmatrix} \begin{pmatrix} \chi'_{j_1}\chi_{j_2}\chi_{j_3} \\ \chi_{j_1}\chi'_{j_2}\chi_{j_3} \\ \chi_{j_1}\chi_{j_2}\chi'_{j_3} \end{pmatrix} \quad (2.35) \\ &= \chi'_{i_1}\chi'_{j_1}\chi_{i_2}\chi_{j_2}\chi_{i_3}\chi_{j_3}D_{11}+ \\ &\quad \chi'_{i_1}\chi_{j_1}(\chi_{i_2}\chi'_{j_2}\chi_{i_3}\chi_{j_3}D_{12} + \chi_{i_2}\chi_{j_2}\chi_{i_3}\chi'_{j_3}D_{13})+ \\ &\quad \chi_{i_1}\chi'_{j_1}(\chi'_{i_2}\chi_{j_2}\chi_{i_3}\chi_{j_3}D_{21} + \chi_{i_2}\chi_{j_2}\chi'_{i_3}\chi_{j_3}D_{31})+ \\ &\quad \chi_{i_1}\chi_{j_1}(\chi'_{i_2}\chi'_{j_2}\chi_{i_3}\chi_{j_3}D_{22} + \chi'_{i_2}\chi_{j_2}\chi_{i_3}\chi'_{j_3}D_{23} + \chi_{i_2}\chi'_{j_2}\chi'_{i_3}\chi_{j_3}D_{32} + \chi_{i_2}\chi_{j_2}\chi'_{i_3}\chi'_{j_3}D_{33}). \end{aligned}$$
$$(2.36)$$

It is now clear that each term in (2.36) is associated to one entry of $\mathcal{D}$. We can thus propose auxiliary functions identified by the indices of each $D_{ab}$. Additionally, after seeing the structure of (2.36) it would be quite useful to have a tool that allows, in a systematic way, to shift between the shape function $\chi_{i_a}$ and its derivative $\chi'_{i_a}$ from one term to another. The inclusion of a binary superscript $\langle s \rangle$ can do such a task,

$$\chi_{i_a}^{\langle s \rangle} = \begin{cases} \chi_{i_a} & \text{if } s = 0, \\ \chi'_{i_a} & \text{if } s = 1. \end{cases} \quad (2.37)$$

Additionally, a really intuitive yet handy function to manage that binary superscript may be a Kronecker delta. For instance, note that any term of (2.36) can be represented as,

$$\chi_{i_1}^{\langle \delta_{1a} \rangle} \chi_{j_1}^{\langle \delta_{1b} \rangle} \chi_{i_2}^{\langle \delta_{2a} \rangle} \chi_{j_2}^{\langle \delta_{2b} \rangle} \chi_{i_3}^{\langle \delta_{3a} \rangle} \chi_{j_3}^{\langle \delta_{3b} \rangle} D_{ab}$$

with $a, b = 1, 2, 3$.

Although this whole article has as an important point of reference the algorithm of sum factorization for the Helmholtz equation (whose stiffness matrix is really similar to a $H^1$ Gram matrix) in the dissertation of Jason Kurtz [18], where each term of the factorized form of the integral conveyed a different auxiliary function, here we decide to provide a single definition of auxiliary functions at each level. These are followed by a sum over $a, b$, hence returning all the terms for the required Gram matrix. This new perspective was applied because one of the goals of the current work is to develop a general framework for the four spaces in the exact sequence. Consider the following definitions,

$$\mathcal{G}^{\text{grad } A}_{ab;i_3j_3}(\xi_1, \xi_2) := \int_0^1 \chi_{i_3}^{\langle \delta_{3a} \rangle}(\xi_3)\chi_{j_3}^{\langle \delta_{3b} \rangle}(\xi_3) D_{ab}(\xi_1, \xi_2, \xi_3)|\mathcal{J}(\xi_1, \xi_2, \xi_3)|d\xi_3, \tag{2.38}$$

$$\mathcal{G}^{\text{grad } B}_{ab;i_2j_2i_3j_3}(\xi_1) := \int_0^1 \chi_{i_2}^{\langle \delta_{2a} \rangle}(\xi_2)\chi_{j_2}^{\langle \delta_{2b} \rangle}(\xi_2) \mathcal{G}^{\text{grad } A}_{ab;i_3j_3}(\xi_1, \xi_2)d\xi_2, \tag{2.39}$$

$$\mathcal{G}^{\text{grad}}_{ab;i_1j_1i_2j_2i_3j_3} := \int_0^1 \chi_{i_1}^{\langle \delta_{1a} \rangle}(\xi_1)\chi_{j_1}^{\langle \delta_{1b} \rangle}(\xi_1) \mathcal{G}^{\text{grad } B}_{ab;i_2j_2i_3j_3}(\xi_1)d\xi_1. \tag{2.40}$$

On the other hand, the first term of the inner product is computed with a slight variation of (2.38)-(2.40), but with no derivative and without multiplying the integrand by $D_{ab}$, therefore that pair of indices is not needed.

$$\bar{\mathcal{G}}^{\text{grad } A}_{i_3j_3}(\xi_1, \xi_2) := \int_0^1 \chi_{i_3}(\xi_3)\chi_{j_3}(\xi_3)|\mathcal{J}(\xi_1, \xi_2, \xi_3)|d\xi_3, \tag{2.41}$$

$$\bar{\mathcal{G}}^{\text{grad } B}_{i_2j_2i_3j_3}(\xi_1) := \int_0^1 \chi_{i_2}(\xi_2)\chi_{j_2}(\xi_2)\bar{\mathcal{G}}^{\text{grad } A}_{i_3j_3}(\xi_1, \xi_2)d\xi_2, \tag{2.42}$$

$$\bar{\mathcal{G}}^{\text{grad}}_{i_1j_1i_2j_2i_3j_3} := \int_0^1 \chi_{i_1}(\xi_1)\chi_{j_1}(\xi_1)\bar{\mathcal{G}}^{\text{grad } B}_{i_2j_2i_3j_3}(\xi_1)d\xi_1. \tag{2.43}$$

The addition of the two parts described above yields,

$$G^{\text{grad}}_{IJ} = \bar{\mathcal{G}}^{\text{grad}}_{i_1j_1i_2j_2i_3j_3} + \sum_{a,b=1}^{3} \mathcal{G}^{\text{grad}}_{ab;i_1j_1i_2j_2i_3j_3}. \tag{2.44}$$

Given the auxiliary functions' definitions for this Gram matrix, some symmetries are possible to be identified so that many extra computations may be avoided.

$$\begin{aligned}
\mathcal{G}^{\text{grad } A}_{ab;i_3j_3} &= \mathcal{G}^{\text{grad } A}_{ba;j_3i_3} & &\text{for all indices,} \\
\mathcal{G}^{\text{grad}}_{ab;i_1j_1i_2j_2i_3j_3} &= \mathcal{G}^{\text{grad}}_{ba;j_1i_1j_2i_2j_3i_3} & &\text{for all indices,} \\
\bar{\mathcal{G}}^{\text{grad } A}_{i_3j_3} &= \bar{\mathcal{G}}^{\text{grad } A}_{j_3i_3} & &\text{for all indices,} \\
\bar{\mathcal{G}}^{\text{grad } B}_{i_2j_2i_3j_3} &= \bar{\mathcal{G}}^{\text{grad } B}_{j_2i_2i_3j_3} &= \bar{\mathcal{G}}^{\text{grad } B}_{i_2j_2j_3i_3} & &\text{for all indices,} \\
\bar{\mathcal{G}}^{\text{grad}}_{i_1j_1i_2j_2i_3j_3} &= \bar{\mathcal{G}}^{\text{grad}}_{j_1i_1i_2j_2i_3j_3} &= \bar{\mathcal{G}}^{\text{grad}}_{i_1j_1j_2i_2i_3j_3} &= \bar{\mathcal{G}}^{\text{grad}}_{i_1j_1i_2j_2j_3i_3} & \text{for all indices.}
\end{aligned}$$

Having in mind the relations above, along with (2.38)-(2.43), we propose Algorithm 5 as the tensorization-based method to compute the $H^1$ Gram matrix, in which a one-dimensional quadrature is applied to every coordinate, as seen in (2.23)-(2.25). In the algorithm, a new hypothetical function named $fillsymmetric3H1$ has the role of filling in the entries that may be obtained just by applying the index symmetries of the introduced auxiliary function arrays $\mathcal{G}^{\text{grad}}$ and $\bar{\mathcal{G}}^{\text{grad } A}$. However, the presence of such a function is not imperative, and the construction of the final matrix could be made inside the algorithm, after appropriate adjustments. This observation applies to the rest of the tensorization-based algorithms presented.

If a uniform order $p$ is assumed in Algorithm 5, we need at least $L = p+1$ for the integration, so the reader may estimate the operation count by summing the number of times the code line that accumulates $\mathcal{G}_{i_1 i j_1 i_2 j_2 i_3 j_3}$ and verify that it is $\mathcal{O}[(p+1)^7]$, two orders of magnitude below Algorithm 4.

---

**Algorithm 5** Tensorization-based computation of the $H^1$ Gram Matrix

---

   **procedure** H1GRAMTENSOR($mid, G^{\text{grad}}$)    ▷ Get $G^{\text{grad}}$ for element No. $mid$ - Tensorization-based algo.
      $p_{\max} \leftarrow \max\{p_1, p_2, p_3\}$
      **call** setquadrature1D($mid, p_{\max}, L, \zeta^l, w^l$)
      $\bar{\mathcal{G}}^{\text{grad}}, \mathcal{G}^{\text{grad}} \leftarrow 0$    ▷ Initialize Gram Matrix
      **for** $l = 1$ to $L$ **do**
         **call** Shape1H1($\zeta^l, p_1, \chi_{i_1}(\zeta^l), \chi'_{i_1}(\zeta^l)\}$)    ▷ Evaluate 1D shape functions at $\zeta^l$
         **for** $j_3 = 0$ to $p_3$ **do**
            **for** $i_3 = j_3$ to $p_3$ **do**
               $\bar{\mathcal{G}}^{\text{grad }B}, \mathcal{G}^{\text{grad }B} \leftarrow 0$
               **for** $m = 1$ to $L$ **do**
                  **call** Shape1H1($\zeta^m, p_2, \chi_{i_2}(\zeta^m), \chi'_{i_2}(\zeta^m)\}$)    ▷ Evaluate 1D shape functions at $\zeta^m$
                  $\bar{\mathcal{G}}^{\text{grad }A}, \mathcal{G}^{\text{grad }A} \leftarrow 0$
                  **for** $n = 1$ to $L$ **do**
                      **call** Shape1H1($\zeta^n, p_3, \chi_{i_3}(\zeta^n), \chi'_{i_3}(\zeta^n)\}$)    ▷ Evaluate 1D shape functions at $\zeta^n$
                      $\boldsymbol{\xi}_{lmn} \leftarrow (\zeta^l, \zeta^m, \zeta^n)$
                      **call** geometry($\boldsymbol{\xi}_{lmn}, mid, \mathbf{x}, \mathcal{J}(\boldsymbol{\xi}_{lmn}), \mathcal{J}^{-1}(\boldsymbol{\xi}_{lmn}), |\mathcal{J}|$)    ▷ Compute jacobian
                      $\mathcal{D} \leftarrow \mathcal{J}^{-1}(\boldsymbol{\xi}_{lmn})\mathcal{J}^{-T}(\boldsymbol{\xi}_{lmn})$
                      $\bar{\mathcal{G}}^{\text{grad }A}_{i_3 j_3} \leftarrow \bar{\mathcal{G}}^{\text{grad }A}_{i_3 j_3} + \chi_{i_3}(\zeta^n)\chi_{j_3}(\zeta^n)|\mathcal{J}|w^n$    ▷ Accumulate to obtain (2.41)
                      **for** $a, b = 1$ to $3$ **do**
                         $\mathcal{G}^{\text{grad }A}_{ab; i_3 j_3} \leftarrow \mathcal{G}^{\text{grad }A}_{ab; i_3 j_3} + \chi_{i_3}^{\langle\delta_{3a}\rangle}(\zeta^n)\chi_{j_3}^{\langle\delta_{3b}\rangle}(\zeta^n)D_{ab}(\zeta^n)|\mathcal{J}|w^n$    ▷ Accumulate to obtain (2.38)
                  **for** $j_2, i_2 = 0$ to $p_2$ **do**
                      $\bar{\mathcal{G}}^{\text{grad }B}_{i_2 j_2 i_3 j_3} \leftarrow \bar{\mathcal{G}}^{\text{grad }B}_{i_2 j_2 i_3 j_3} + \chi_{i_2}(\zeta^m)\chi_{j_2}(\zeta^m)\hat{\mathcal{G}}^{\text{grad }A}_{i_3 j_3}(\zeta^l, \zeta^m)w^m$    ▷ By (2.42)
                      **for** $a, b = 1$ to $3$ **do**
                         $\mathcal{G}^{\text{grad }B}_{ab; i_2 j_2 i_3 j_3} \leftarrow \mathcal{G}^{\text{grad }B}_{ab; i_2 j_2 i_3 j_3} + \chi_{i_2}^{\langle\delta_{2a}\rangle}(\zeta^m)\chi_{j_2}^{\langle\delta_{2b}\rangle}(\zeta^m)\mathcal{G}^{\text{grad }A}_{ab; i_3 j_3}(\zeta^l, \zeta^m)w^m$    ▷ By (2.39)
            **for** $j_2, i_2 = 0$ to $p_2$ **do**
                **for** $j_1, i_1 = 0$ to $p_1$ **do**
                    $I = i_1 + (p_1+1)i_2 + (p_1+1)(p_2+1)i_3$
                    $J = j_1 + (p_1+1)j_2 + (p_1+1)(p_2+1)j_3$
                    **if** $J \geq I$ **then**
                        $\bar{\mathcal{G}}^{\text{grad}}_{i_1 j_1 i_2 j_2 i_3 j_3} \leftarrow \bar{\mathcal{G}}^{\text{grad}}_{i_1 j_1 i_2 j_2 i_3 j_3} + \chi_{i_1}(\zeta^l)\chi_{j_1}(\zeta^l)\bar{\mathcal{G}}^{\text{grad }B}_{i_2 j_2 i_3 j_3}(\zeta^l)w^l$    ▷ (2.43)
                        **for** $a, b = 1$ to $3$ **do**
                            $\mathcal{G}^{\text{grad}}_{ab; i_1 j_1 i_2 j_2 i_3 j_3} \leftarrow \mathcal{G}^{\text{grad}}_{ab; i_1 j_1 i_2 j_2 i_3 j_3} + \chi_{i_1}^{\langle\delta_{1a}\rangle}(\zeta^l)\chi_{j_1}^{\langle\delta_{1b}\rangle}(\zeta^l)\mathcal{G}^{\text{grad }B}_{ab; i_2 j_2 i_3 j_3}(\zeta^l)w^l$    ▷ (2.40)
      **call** fillsymmetric3H1($p_1, p_2, p_3, \mathcal{G}^{\text{grad}}, \bar{\mathcal{G}}^{\text{grad}}, G^{\text{grad}}$)    ▷ Adds terms and fills matrix using symmetries
      **return** $G^{\text{grad}}$

---

Exactly as in the $L^2$ case, matrices $\bar{\mathcal{G}}^{\text{grad}}_{i_1 j_1 i_2 j_2 i_3 j_3}$ and $\mathcal{G}^{\text{grad}}_{ab; i_1 j_1 i_2 j_2 i_3 j_3}$ need not be constructed and stored because that would demand too much memory. To actually translate these lines into the code, for each combination of $i_1 j_1 i_2 j_2 i_3 j_3$ the values for such arrays can be accumulated into temporary variables and then added into the corresponding spot of $G^{\text{grad}}$. This same idea holds when implementing the algorithms for the next two energy spaces.

## 2.5 Space $H(\text{div})$

This third space contains functions whose divergence is square integrable, i.e,

$$H(\text{div}, \mathcal{K}) = \left\{ V \in \boldsymbol{L}^2(\mathcal{K}) \; : \text{div}\, V \in L^2(\mathcal{K}) \right\}. \tag{2.45}$$

The corresponding inner product is:

$$(V, W)_{H(\text{div}, \mathcal{K})} := (V, W)_\mathcal{K} + (\text{div}\, V, \text{div}\, W)_\mathcal{K} \quad \forall V, W \in H(\text{div}, \mathcal{K}); \tag{2.46}$$

where the spatial-coordinates divergence is:

$$\text{div}\, V = \partial_1 V_1 + \partial_2 V_2 + \partial_3. \tag{2.47}$$

Furthermore, the divergence in the master element coordinates of a function $\hat{W}$ is:

$$\widehat{\text{div}}\hat{W} = \hat{\partial}_1 \hat{W}_1 + \hat{\partial}_2 \hat{W}_2 + \hat{\partial}_3 \hat{W}_3. \tag{2.48}$$

In a similar way to the previous problem, we have that if $V \in H(\text{div}, \mathcal{K})$, then $\text{div}\, V \in L^2(\mathcal{K})$; and that if $V = T^{\text{div}} \hat{V}$ for some $\hat{V} \in H(\text{div}, \hat{\mathcal{K}})$, then $\text{div}\, V = T\widehat{\text{div}}\hat{V}$. Consequently, $\widehat{\text{div}}\hat{V} \in L^2(\hat{\mathcal{K}})$.

Now, let the order of the shape functions for the master hexahedron be $(p_1, p_2, p_3)$, in the sense of the exact sequence. Consider a basis for $V^p$, $\{\vartheta_I\}_{I=0}^{\dim V^p - 1}$, where $\dim V^p = (p_1 + 1)p_2 p_3 + (p_2 + 1)p_3 p_1 + (p_3 + 1)p_1 p_2$. Central to the understanding of this third space is to remind that its elements are vector-valued functions. As in the $H^1$ case, the definitions of both $T^{\text{div}}$ and $T$ need to be invoked in (2.49). With this in mind, we take the same steps as above to achieve an expression suitable for tensorization for computing the $H(\text{div})$ Gram matrix, below denoted $G^{\text{div}}$. Thus, for any pair of integers $0 \leq I, J < \dim V^p$ the inner product is obtained as follows:

$$\begin{aligned}
G^{\text{div}}_{IJ} &= (\vartheta_I, \vartheta_J)_{H(\text{div}, \mathcal{K})} \\
&= \int_\mathcal{K} \vartheta_I(\boldsymbol{x})^\mathsf{T} \vartheta_J(\boldsymbol{x}) d^3 \boldsymbol{x} + \int_\mathcal{K} \text{div}\, \vartheta_I(\boldsymbol{x}) \, \text{div}\, \vartheta_J(\boldsymbol{x}) d^3 \boldsymbol{x} \\
&= \int_{\hat{\mathcal{K}}} [\vartheta_I \circ \mathbf{x}_\mathcal{K}(\boldsymbol{\xi})]^\mathsf{T} [\vartheta_J \circ \mathbf{x}_\mathcal{K}(\boldsymbol{\xi})] \, |\mathcal{J}(\boldsymbol{\xi})| d^3 \boldsymbol{\xi} + \int_\mathcal{K} \text{div}\, \vartheta_I \circ \mathbf{x}_\mathcal{K}(\boldsymbol{\xi}) \, \text{div}\, \vartheta_J \circ \mathbf{x}_\mathcal{K}(\boldsymbol{\xi}) |\mathcal{J}(\boldsymbol{\xi})| d^3 \boldsymbol{\xi} \\
&= \int_{\hat{\mathcal{K}}} \left[ T^{\text{div}} \hat{\vartheta}_I \circ \mathbf{x}_\mathcal{K}(\boldsymbol{\xi}) \right]^\mathsf{T} \left[ T^{\text{div}} \hat{\vartheta}_J \circ \mathbf{x}_\mathcal{K}(\boldsymbol{\xi}) \right] |\mathcal{J}(\boldsymbol{\xi})| d^3 \boldsymbol{\xi} + \\
&\quad \int_\mathcal{K} \left[ \text{div}\left( T^{\text{div}} \hat{\vartheta}_I \right) \circ \mathbf{x}_\mathcal{K}(\boldsymbol{\xi}) \right] \left[ \text{div}\left( T^{\text{div}} \hat{\vartheta}_J \right) \circ \mathbf{x}_\mathcal{K}(\boldsymbol{\xi}) \right] |\mathcal{J}(\boldsymbol{\xi})| d^3 \boldsymbol{\xi}
\end{aligned}$$

$$
\begin{aligned}
&= \int_{\hat{\mathcal{K}}} \left[\left(|\mathcal{J}|^{-1}\mathcal{J}\hat{\vartheta}_I \circ \mathbf{x}_{\mathcal{K}}^{-1}\right) \circ \mathbf{x}_{\mathcal{K}}(\boldsymbol{\xi})\right]^{\mathsf{T}} \left[\left(|\mathcal{J}|^{-1}\mathcal{J}\hat{\vartheta}_J \circ \mathbf{x}_{\mathcal{K}}^{-1}\right) \circ \mathbf{x}_{\mathcal{K}}(\boldsymbol{\xi})\right] |\mathcal{J}(\boldsymbol{\xi})| d^3\boldsymbol{\xi} + \\
&\quad \int_{\mathcal{K}} \left[\left(T\widehat{\mathrm{div}}\hat{\vartheta}_I\right) \circ \mathbf{x}_{\mathcal{K}}(\boldsymbol{\xi})\right] \left[\left(T\widehat{\mathrm{div}}\hat{\vartheta}_J\right) \circ \mathbf{x}_{\mathcal{K}}(\boldsymbol{\xi})\right] |\mathcal{J}(\boldsymbol{\xi})| d^3\boldsymbol{\xi} \\
&= \int_{\hat{\mathcal{K}}} \hat{\vartheta}_I(\boldsymbol{\xi})^{\mathsf{T}} \mathcal{C}(\boldsymbol{\xi}) \hat{\vartheta}_J(\boldsymbol{\xi}) |\mathcal{J}^{-1}(\boldsymbol{\xi})| d^3\boldsymbol{\xi} + \\
&\quad \int_{\mathcal{K}} \left[\left((|\mathcal{J}|^{-1}\widehat{\mathrm{div}}\hat{\vartheta}_I) \circ \mathbf{x}_{\mathcal{K}}^{-1}\right) \circ \mathbf{x}_{\mathcal{K}}(\boldsymbol{\xi})\right] \left[\left((|\mathcal{J}|^{-1}\widehat{\mathrm{div}}\hat{\vartheta}_J) \circ \mathbf{x}_{\mathcal{K}}^{-1}\right) \circ \mathbf{x}_{\mathcal{K}}(\boldsymbol{\xi})\right] |\mathcal{J}(\boldsymbol{\xi})| d^3\boldsymbol{\xi} \\
&= \int_{\hat{\mathcal{K}}} \hat{\vartheta}_I(\boldsymbol{\xi})^{\mathsf{T}} \mathcal{C}(\boldsymbol{\xi}) \hat{\vartheta}_J(\boldsymbol{\xi}) |\mathcal{J}(\boldsymbol{\xi})|^{-1} d^3\boldsymbol{\xi} + \int_{\mathcal{K}} \widehat{\mathrm{div}}\hat{\vartheta}_I(\boldsymbol{\xi}) \widehat{\mathrm{div}}\hat{\vartheta}_J(\boldsymbol{\xi}) |\mathcal{J}(\boldsymbol{\xi})|^{-1} d^3\boldsymbol{\xi} \\
&= \int_0^1 \int_0^1 \int_0^1 \left\{\hat{\vartheta}_I(\boldsymbol{\xi})^{\mathsf{T}} \mathcal{C}(\boldsymbol{\xi}) \hat{\vartheta}_J(\boldsymbol{\xi}) + \widehat{\mathrm{div}}\hat{\vartheta}_I(\boldsymbol{\xi}) \widehat{\mathrm{div}}\hat{\vartheta}_J(\boldsymbol{\xi})\right\} |\mathcal{J}(\boldsymbol{\xi})|^{-1} d\xi_3 d\xi_2 d\xi_1.
\end{aligned} \quad (2.49)
$$

where $\mathcal{C} := \mathcal{J}^{\mathsf{T}}\mathcal{J} = (\mathcal{J}^{-1}\mathcal{J}^{-T})^{-1} = \mathcal{D}^{-1}$ contains the following entries:

$$
\mathcal{C} = \begin{pmatrix} C_{11} & C_{12} & C_{13} \\ C_{21} & C_{22} & C_{23} \\ C_{31} & C_{32} & C_{33} \end{pmatrix}, \quad (2.50)
$$

with $C_{12} = C_{21}$, $C_{13} = C_{31}$, $C_{32} = C_{23}$, and every entry depends on $\boldsymbol{\xi}$.

The conventional algorithm for (2.49) is presented as Algorithm 6.

---

**Algorithm 6** Conventional computation of the $H(\mathrm{div})$ Gram Matrix
---

    **procedure** HDIVGRAM($mid, G^{\mathrm{div}}$) ▷ Compute matrix $G^{\mathrm{div}}$ for element No. $mid$ - Conventional algorithm
        **call** setquadrature3D($mid, p_1, p_2, p_3, L, \boldsymbol{\xi}_{lmn}, w_{lmn}$)
        $G^{\mathrm{div}} \leftarrow 0$         ▷ Initialize Gram Matrix
        **for** $l, m, n = 1$ to $L$ **do**
            **call** Shape3Hdiv($\boldsymbol{\xi}_{lmn}, p_1, p_2, p_3, \{\hat{\vartheta}_I(\boldsymbol{\xi}_{lmn})\}, \{\widehat{\mathrm{div}}\hat{\vartheta}_I(\boldsymbol{\xi}_{lmn})\}$)   ▷ 3D shape functions at $\boldsymbol{\xi}_{lmn}$
            **call** geometry($\boldsymbol{\xi}_{lmn}, mid, \mathbf{x}, \mathcal{J}(\boldsymbol{\xi}_{lmn}), \mathcal{J}^{-1}(\boldsymbol{\xi}_{lmn}), |\mathcal{J}|$)   ▷ Compute jacobian
            $\mathcal{C} \leftarrow \mathcal{J}^{\mathsf{T}}(\boldsymbol{\xi}_{lmn}) \mathcal{J}(\boldsymbol{\xi}_{lmn})$
            **for** $J = 0$ to $\dim V^p - 1$ **do**
                **for** $I = J$ to $\dim V^p - 1$ **do**
                    $G^{\mathrm{div}}_{IJ} \leftarrow G^{\mathrm{div}}_{IJ} + \left\{\hat{\vartheta}_I(\boldsymbol{\xi}_{lmn})^{\mathsf{T}} \mathcal{C} \hat{\vartheta}_J(\boldsymbol{\xi}_{lmn}) + \widehat{\mathrm{div}}\hat{\vartheta}_I(\boldsymbol{\xi}_{lmn}) \widehat{\mathrm{div}}\hat{\vartheta}_J(\boldsymbol{\xi}_{lmn})\right\} |\mathcal{J}|^{-1} w_{lmn}$
    **return** $G^{\mathrm{div}}$

---

The tensorization-based integration for the present situation begins by defining a tensor-product shape function. As $\hat{\vartheta}_I, \hat{\vartheta}_I \in \hat{V}^p = \mathcal{Q}^{p_1, p_2-1, p_3-1}(\mathcal{I}^3) \times \mathcal{Q}^{p_1-1, p_2, p_3-1}(\mathcal{I}^3) \times \mathcal{Q}^{p_1-1, p_2-1, p_3}(\mathcal{I}^3)$, it follows that:

$$
\begin{aligned}
\hat{\vartheta}_I(\xi_1, \xi_2, \xi_3) &= \mu_{a,1;i_1}(\xi_1) \mu_{a,2;i_2}(\xi_2) \mu_{a,3;i_3}(\xi_3) \hat{\mathbf{e}}_a \\
\hat{\vartheta}_J(\xi_1, \xi_2, \xi_3) &= \mu_{b,1;j_1}(\xi_1) \mu_{b,2;j_2}(\xi_2) \mu_{b,3;j_3}(\xi_3) \hat{\mathbf{e}}_b.
\end{aligned} \quad (2.51)
$$

where $\hat{\mathbf{e}}_a$, $a = 1, 2, 3$, is the $a$-th canonical Cartesian unit vector in the master space; the univariate polyno-mials $\mu_{a,d;i_d}(\xi_d)$, $d = 1, 2, 3$, are determined through the rule

$$\mu_{a,d;i_d} = \begin{cases} \chi_{i_d} & \text{if } a = d, \\ \nu_{i_d} & \text{otherwise.} \end{cases} \qquad (2.52)$$

where it is recalled that $\{\chi_{i_a}\}_{i_a=0}^{p_a}$ are a hierarchical basis of $W_\mathcal{I}^{p_a}$, while $\{\nu_{i_a}\}_{i_a=0}^{p_a-1}$ form a basis for $Y_\mathcal{I}^{p_a}$ (for $a = 1, 2, 3$). In this new scenario, the integer indices $i_a, j_a$ depend on the vector component where $\hat{\vartheta}_I$ lies. The correspondence formula between the tensor-product index $I$ and the one-dimensional ones could be of the form

$$I = \begin{cases} i_1 + (p_1 + 1)i_2 + (p_1 + 1)p_2 i_3 & \text{if a=1}, \\ (p_1 + 1)p_2 p_3 + i_1 + p_1 i_2 + p_1(p_2 + 1)i_3 & \text{if a=2}, \\ (p_1 + 1)p_2 p_3 + p_1(p_2 + 1)p_3 + i_1 + p_1 i_2 + p_1 p_2 i_3 & \text{if a=3}. \end{cases} \qquad (2.53)$$

Furthermore, it should be acknowledged that if $\{\chi_{i_a}\}_{i_a=0}^{p_a}$ is a hierarchical basis of $\mathcal{P}^{p_a}(\mathcal{I})$, then $\{\chi'_{i_a}\}_{i_a=1}^{p_a}$ is a hierarchical basis of $\mathcal{P}^{p_a-1}(\mathcal{I}) = Y_\mathcal{I}^{p_a}$. We can therefore use this basis by letting

$$\nu_{i_d} = \chi'_{i_d+1} \text{ for } i_d = 0, ..., p_d - 1, \qquad (2.54)$$

having $p_d = p_1, p_2, p_3$.

By combining (2.51), (2.52) and (2.54) the tensor-product $H(\text{div})$ shape functions and its master-coordinate divergence can be rewritten as:

$$\hat{\vartheta}_I(\xi_1, \xi_2, \xi_3) = \chi_{i_a}(\xi_a)\chi'_{i_{a+1}+1}(\xi_{a+1})\chi'_{i_{a+2}+1}(\xi_{a+2})\hat{\mathbf{e}}_a, \qquad (2.55)$$

$$\widehat{\text{div}}\hat{\vartheta}_I(\xi_1, \xi_2, \xi_3) = \chi'_{i_a}(\xi_a)\chi'_{i_{a+1}+1}(\xi_{a+1})\chi'_{i_{a+2}+1}(\xi_{a+2}), \qquad (2.56)$$

with $0 \le i_a \le p_a$, $0 \le i_{a+1} < p_{a+1}$, $0 \le i_{a+2} < p_{a+2}$; where $a = 1, 2, 3$ and the indices $(a, a+1, a+2)$ are interpreted as a cyclic permutation of $(1, 2, 3)$.

Notice that if we introduce an additional identifier to the index, we can restate the shape function in terms of $(\xi_1, \xi_2, \xi_3)$ which is desirable before the integral nesting of the tensorization-based algorithm. As above, a Kronecker delta associated to the vector component where the shape function lies can do that work for us. Additionally, the notation convention to know if we need $\chi$ or its derivative turns again useful. Making use of (2.37) jointly with the Kronecker delta idea for the index, we can obtain new equations for $\hat{\vartheta}_I$ and $\widehat{\text{div}}\hat{\vartheta}_I$, for which the effect of (the component index) $a$ on each univariate polynomial is obtained more systematically than in (2.55) and (2.56), respectively.

$$\hat{\vartheta}_I(\xi_1, \xi_2, \xi_3) = \chi_{i_1+1-\delta_{1a}}^{\langle 1-\delta_{1a}\rangle}(\xi_1)\chi_{i_2+1-\delta_{2a}}^{\langle 1-\delta_{2a}\rangle}(\xi_2)\chi_{i_3+1-\delta_{3a}}^{\langle 1-\delta_{3a}\rangle}(\xi_3)\hat{\mathbf{e}}_a, \qquad (2.57)$$

$$\widehat{\text{div}}\hat{\vartheta}_I(\xi_1, \xi_2, \xi_3) = \chi'_{i_1+1-\delta_{1a}}(\xi_1)\chi'_{i_2+1-\delta_{2a}}(\xi_2)\chi'_{i_3+1-\delta_{3a}}(\xi_3). \qquad (2.58)$$

Using (2.58) in the second term of the Gram matrix entry derived in (2.49) we get to

$$\int_{\hat{\mathcal{K}}} \widehat{\text{div}}\hat{\vartheta}_I \widehat{\text{div}}\hat{\vartheta}_J |\mathcal{J}|^{-1} d^3 \boldsymbol{\xi} = \int_{\hat{\mathcal{K}}} \chi'_{i_1+1-\delta_{1a}}\chi'_{i_2+1-\delta_{2a}}\chi'_{i_3+1-\delta_{3a}}\chi'_{j_1+1-\delta_{1b}}\chi'_{j_2+1-\delta_{2b}}\chi'_{j_3+1-\delta_{3b}}|\mathcal{J}|^{-1} d^3 \boldsymbol{\xi}, \qquad (2.59)$$

which, due to (2.54), is indeed an $L^2$ inner product just as in (2.19). Notice that in (2.59) we abolished the writing of arguments, in order to shorten the number of symbols present, but it clearly holds that each univariate shape function is dependent on the $\boldsymbol{\xi}$ component with its same index. As we have presented an algorithm for this type of integral, let us focus on the first term of (2.49) instead, and we'll later incorporate both of them in a single algorithm.

Take (2.57) as a representation of the $H(\text{div})$ shape functions $\hat{\vartheta}_I, \hat{\vartheta}_J$. Insert this representation into the first term of (2.49), and rearrange in 1D integrals to obtain a nested integral,

$$\int_{\hat{\mathcal{K}}} \hat{\vartheta}_I^\mathsf{T} \mathcal{C} \hat{\vartheta}_J |\mathcal{J}|^{-1} d^3\boldsymbol{\xi} = \int_{\hat{\mathcal{K}}} \chi_{i_1+1-\delta_{1a}}^{\langle 1-\delta_{1a}\rangle} \chi_{i_2+1-\delta_{2a}}^{\langle 1-\delta_{2a}\rangle} \chi_{i_3+1-\delta_{3a}}^{\langle 1-\delta_{3a}\rangle} \hat{\mathbf{e}}_a^\mathsf{T} \mathcal{C} \chi_{j_1+1-\delta_{1b}}^{\langle 1-\delta_{1b}\rangle} \chi_{j_2+1-\delta_{2b}}^{\langle 1-\delta_{2b}\rangle} \chi_{j_3+1-\delta_{3b}}^{\langle 1-\delta_{3b}\rangle} \hat{\mathbf{e}}_b |\mathcal{J}|^{-1} d^3\boldsymbol{\xi}$$

$$= \int_0^1 \chi_{i_1+1-\delta_{1a}}^{\langle 1-\delta_{1a}\rangle} \chi_{j_1+1-\delta_{1b}}^{\langle 1-\delta_{1b}\rangle} \left\{ \int_0^1 \chi_{i_2+1-\delta_{2a}}^{\langle 1-\delta_{2a}\rangle} \chi_{j_2+1-\delta_{2b}}^{\langle 1-\delta_{2b}\rangle} \left[ \int_0^1 \chi_{i_3+1-\delta_{3a}}^{\langle 1-\delta_{3a}\rangle} \chi_{j_3+1-\delta_{3b}}^{\langle 1-\delta_{3b}\rangle} C_{ab} |\mathcal{J}|^{-1} d\xi_3 \right] d\xi_2 \right\} d\xi_1, \quad (2.60)$$

where we have used the fact that $\hat{\mathbf{e}}_a^\mathsf{T} \mathcal{C} \hat{\mathbf{e}}_b = C_{ab}$.

This is where the corresponding auxiliary functions for $H(\text{div})$ are introduced.

$$\mathcal{G}_{ab;i_3j_3}^{\text{div}\,A}(\xi_1, \xi_2) := \int_0^1 \chi_{i_3+1-\delta_{3a}}^{\langle 1-\delta_{3a}\rangle}(\xi_3) \chi_{j_3+1-\delta_{3b}}^{\langle 1-\delta_{3b}\rangle}(\xi_3) C_{ab}(\xi_1, \xi_2, \xi_3) |\mathcal{J}(\xi_1, \xi_2, \xi_3)|^{-1} d\xi_3, \quad (2.61)$$

$$\mathcal{G}_{ab;i_2j_2i_3j_3}^{\text{div}\,B}(\xi_1) := \int_0^1 \chi_{i_2+1-\delta_{2a}}^{\langle 1-\delta_{2a}\rangle}(\xi_2) \chi_{j_2+1-\delta_{2b}}^{\langle 1-\delta_{2b}\rangle}(\xi_2) \mathcal{G}_{ab;i_3j_3}^{\text{div}\,A}(\xi_1, \xi_2) d\xi_2, \quad (2.62)$$

$$\mathcal{G}_{ab;i_1j_2i_2j_2i_3j_3}^{\text{div}}(\xi_1) := \int_0^1 \chi_{i_1+1-\delta_{1a}}^{\langle 1-\delta_{1a}\rangle}(\xi_1) \chi_{j_1+1-\delta_{1b}}^{\langle 1-\delta_{1b}\rangle}(\xi_1) \mathcal{G}_{ab;i_2j_2i_3j_3}^{\text{div}\,B}(\xi_1) d\xi_1. \quad (2.63)$$

In these newly defined auxiliary function arrays for the first term of the $H(\text{div})$ Gram matrix, some symmetries arise as well as they did in the $L^2$ and $H^1$ cases. We don't explicitly state them now, but the reader is encouraged to explore such properties, having as a main criterion the fact that the Gram matrix as a whole is symmetric.

Recall that the second term of the inner product, written in tensor-product form in (2.59), must be added to (2.63). The calculation of that term is practically a copy of Algorithm 2, but some caution has to be taken with the indices. We'll denote those slightly modified auxiliary arrays $\widetilde{\mathcal{G}}^{\text{div}\,A}, \widetilde{\mathcal{G}}^{\text{div}\,B}, \widetilde{\mathcal{G}}^{\text{div}}$. Next, by implementing a 1D quadrature for the approximation of each auxiliary function (2.61)-(2.63), similar to those in (2.23)-(2.25), we get the third tensorization-based algorithm for the computation of a hexahedral element's Gram matrix, Algorithm 7. This algorithm takes into account the major symmetry of $G^{\text{div}}$, but if it is wanted to use some symmetries of the auxiliary functions, those would be implemented within the hypothetical function $fillsymmetric3Hdiv$ (which is not imperative, as explained above). The sum of $\mathcal{G}^{\text{div}}$ and $\widetilde{\mathcal{G}}^{\text{div}}$ could also be made within such a function.

## 2.6 Space $H(\text{curl})$

The last energy space needs its functions to have a square-integrable curl. It is defined as

$$H(\text{curl}, \mathcal{K}) = \left\{ E \in \boldsymbol{L}^2(\mathcal{K}) \ : \text{curl}\, E \in \boldsymbol{L}^2(\mathcal{K}) \right\}. \quad (2.64)$$

**Algorithm 7** Tensorization-based computation of the $H(\text{div})$ Gram Matrix

**procedure** HDIVGRAMTENSOR($mid, G^{\text{div}}$) ▷ $G^{\text{div}}$ for element No. $mid$ -Tensorization - based algorithm
    $p_{\max} \leftarrow \max\{p_1, p_2, p_3\}$
    **call** setquadrature1D($mid, p_{\max}, L, \zeta^l, w^l$)
    $\mathcal{G}^{\text{div}}, \widetilde{\mathcal{G}}^{\text{div}} \leftarrow 0$     ▷ Initialize Gram Matrix
    **for** $l = 1$ to $L$ **do**
        **call** Shape1H1($\zeta^l, p_1, \chi_{i_1}(\zeta^l), \chi'_{i_1}(\zeta^l)\}$)     ▷ Evaluate 1D shape functions at $\zeta^l$
        **for** $j_3 = 0$ to $p_3$ **do**
            **for** $i_3 = j_3$ to $p_3$ **do**
                $\mathcal{G}^{\text{div }B}, \widetilde{\mathcal{G}}^{\text{div }B} \leftarrow 0$
                **for** $m = 1$ to $L$ **do**
                      **call** Shape1H1($\zeta^m, p_2, \chi_{i_2}(\zeta^m), \chi'_{i_2}(\zeta^m)\}$)     ▷ Evaluate 1D shape functions at $\zeta^m$
                      $\mathcal{G}^{\text{div }A}, \widetilde{\mathcal{G}}^{\text{div }A} \leftarrow 0$
                      **for** $n = 1$ to $L$ **do**
                          **call** Shape1H1($\zeta^n, p_3, \chi_{i_3}(\zeta^n), \chi'_{i_3}(\zeta^n)\}$)     ▷ Evaluate 1D shape functions at $\zeta^n$
                          $\boldsymbol{\xi}_{lmn} \leftarrow (\zeta^l, \zeta^m, \zeta^n)$
                          **call** geometry($\boldsymbol{\xi}_{lmn}, mid, \mathbf{x}, \mathcal{J}(\boldsymbol{\xi}_{lmn}), \mathcal{J}^{-1}(\boldsymbol{\xi}_{lmn}), |\mathcal{J}|$)     ▷ Compute jacobian
                          $\mathcal{C} \leftarrow \mathcal{J}^\mathsf{T}(\boldsymbol{\xi}_{lmn})\mathcal{J}(\boldsymbol{\xi}_{lmn})$
                          **for** $a, b = 1$ to $3$ **do**
                              **if** $j_3 + 1 - \delta_{3b} \leq p_3$ and $i_3 + 1 - \delta_{3a} \leq p_3$ **then**     ▷ Avoids extra computations
                                  $\mathcal{G}^{\text{div }A}_{ab;i_3j_3} \leftarrow \mathcal{G}^{\text{div }A}_{ab;i_3j_3} + \chi^{\langle 1-\delta_{3a}\rangle}_{i_3+1-\delta_{3a}}(\zeta^n)\chi^{\langle 1-\delta_{3b}\rangle}_{j_3+1-\delta_{3b}}(\zeta^n)C_{ab}|\mathcal{J}|^{-1}w^n$     ▷ To obtain (2.61)
                                  $\widetilde{\mathcal{G}}^{\text{div }A}_{ab;i_3j_3} \leftarrow \widetilde{\mathcal{G}}^{\text{div }A}_{ab;i_3j_3} + \chi'_{i_3+1-\delta_{3a}}(\zeta^n)\chi'_{j_3+1-\delta_{3b}}(\zeta^n)|\mathcal{J}|^{-1}w^n$
                      **for** $j_2, i_2 = 0$ to $p_2$ **do**
                          **for** $a, b = 1$ to $3$ **do**
                              **if** $j_2 + 1 - \delta_{2b} \leq p_2$ and $i_2 + 1 - \delta_{2a} \leq p_2$ **then**     ▷ Avoids extra computations
                                  $\mathcal{G}^{\text{div }B}_{ab;i_2j_2i_3j_3} \leftarrow \mathcal{G}^{\text{div }B}_{ab;i_2j_2i_3j_3} + \chi^{\langle 1-\delta_{2a}\rangle}_{i_2+1-\delta_{2a}}(\zeta^m)\chi^{\langle 1-\delta_{2b}\rangle}_{j_2+1-\delta_{2b}}(\zeta^m)\mathcal{G}^{\text{div }A}_{ab;i_3j_3}w^m$     ▷ (2.62)
                                  $\widetilde{\mathcal{G}}^{\text{div }B}_{ab;i_2j_2i_3j_3} \leftarrow \widetilde{\mathcal{G}}^{\text{div }B}_{ab;i_2j_2i_3j_3} + \chi'_{i_2+1-\delta_{2a}}(\zeta^m)\chi'_{j_2+1-\delta_{2b}}(\zeta^m)\widetilde{\mathcal{G}}^{\text{div }A}_{ab;i_3j_3}w^m$
            **for** $j_2, i_2 = 0$ to $p_2$ **do**
                **for** $j_1, i_1 = 0$ to $p_1$ **do**
                      Determine $I, J$ through (2.53)
                      **if** $J \geq I$ **then**
                          **for** $a, b = 1$ to $3$ **do**
                              **if** $j_1 + 1 - \delta_{1b} \leq p_1$ and $i_1 + 1 - \delta_{1a} \leq p_1$ **then**     ▷ Avoids extra computations
                                $\mathcal{G}^{\text{div}}_{ab;i_1j_1i_2j_2i_3j_3} \leftarrow \mathcal{G}^{\text{div}}_{ab;i_1j_1i_2j_2i_3j_3} + \chi^{\langle 1-\delta_{1a}\rangle}_{i_1+1-\delta_{1a}}(\zeta^l)\chi^{\langle 1-\delta_{1b}\rangle}_{j_1+1-\delta_{1b}}(\zeta^l)\mathcal{G}^{\text{div }B}_{ab;i_2j_2i_3j_3}w^l$
                                  $\widetilde{\mathcal{G}}^{\text{div}}_{ab;i_1j_1i_2j_2i_3j_3} \leftarrow \widetilde{\mathcal{G}}^{\text{div}}_{ab;i_1j_1i_2j_2i_3j_3} + \chi'_{i_1+1-\delta_{1a}}(\zeta^l)\chi'_{j_1+1-\delta_{1b}}(\zeta^l)\widetilde{\mathcal{G}}^{B}_{ab;i_2j_2i_3j_3}w^l$
    **call** fillsymmetric3Hdiv($p_1, p_2, p_3, \mathcal{G}^{\text{div}}, \widetilde{\mathcal{G}}^{\text{div}}, G^{\text{div}}$)     ▷ Adds terms and fills matrix using symmetries
    **return** $G^{\text{div}}$

The inner product for this space is

$$(E, F)_{H(\text{curl}, \mathcal{K})} := (E, F)_\mathcal{K} + (\text{curl } E, \text{curl } F)_\mathcal{K} \quad \forall E, F \in H(\text{curl}, \mathcal{K}). \tag{2.65}$$

In this definition the curl is given by:

$$\operatorname{curl} E = (\partial_2 E_3 - \partial_3 E_2, \partial_3 E_1 - \partial_1 E_3, \partial_1 E_2 - \partial_2 E_1). \tag{2.66}$$

Its counterpart in master coordinates is:

$$\widehat{\operatorname{curl}}\hat{F} = \left(\hat{\partial}_2 \hat{F}_3 - \hat{\partial}_3 \hat{F}_2, \hat{\partial}_3 \hat{F}_1 - \hat{\partial}_1 \hat{F}_3, \hat{\partial}_1 \hat{F}_2 - \hat{\partial}_2 \hat{F}_1\right). \tag{2.67}$$

In a similar way to the previous problem, we have that if $E \in H(\operatorname{curl}, \mathcal{K})$, then $\operatorname{curl} E \in H(\operatorname{div}, \mathcal{K})$ and that if $E = T^{\operatorname{curl}}\hat{E}$ for some $\hat{E} \in H(\operatorname{curl}, \hat{\mathcal{K}})$, then $\operatorname{curl} E = T^{\operatorname{div}}\widehat{\operatorname{curl}}\hat{E}$. Consequently, $\widehat{\operatorname{curl}}\hat{E} \in H(\operatorname{div}, \hat{\mathcal{K}})$.

Now, let the order of the shape functions for the master hexahedron be $(p_1, p_2, p_3)$, in the sense of the exact sequence. Consider a basis for $Q^p$, $\{\psi_I\}_{I=0}^{\dim Q^p - 1}$, where $\dim Q^p = p_1(p_2+1)(p_3+1) + (p_1+1)p_2(p_3+1) + (p_1+1)(p_2+1)p_3$. In this fourth energy space the elements are vector-valued functions too. As in the $H^1$ and $H(\operatorname{div})$ cases, the definitions of both $T^{\operatorname{curl}}$ and $T^{\operatorname{div}}$ need to be invoked during the derivation of (2.68). Taking in consideration the previous problems, we take the same steps in order to obtain an expression suitable for the tensorization-based integration of the $H(\operatorname{curl})$ Gram matrix, hereinafter called $G^{\operatorname{curl}}$. Thus, for any pair of integers $0 \le I, J < \dim Q^p$ the inner product is equivalent to the following expressions,

$$\begin{aligned}
G^{\operatorname{curl}}_{IJ} &= (\psi_I, \psi_J)_{H(\operatorname{curl},\mathcal{K})} \\
&= \int_{\mathcal{K}} \psi_I(\boldsymbol{x})^{\mathsf{T}} \psi_J(\boldsymbol{x}) d^3\boldsymbol{x} + \int_{\mathcal{K}} [\operatorname{curl} \psi_I(\boldsymbol{x})]^{\mathsf{T}} [\operatorname{curl} \psi_J(\boldsymbol{x})] d^3\boldsymbol{x} \\
&= \int_{\hat{\mathcal{K}}} [\psi_I \circ \mathbf{x}_{\mathcal{K}}(\boldsymbol{\xi})]^{\mathsf{T}} [\psi_J \circ \mathbf{x}_{\mathcal{K}}(\boldsymbol{\xi})] |\mathcal{J}(\boldsymbol{\xi})| d^3\boldsymbol{\xi} + \int_{\mathcal{K}} [\operatorname{curl} \psi_I \circ \mathbf{x}_{\mathcal{K}}(\boldsymbol{\xi})]^{\mathsf{T}} [\operatorname{curl} \psi_J \circ \mathbf{x}_{\mathcal{K}}(\boldsymbol{\xi})] |\mathcal{J}(\boldsymbol{\xi})| d^3\boldsymbol{\xi} \\
&= \int_{\hat{\mathcal{K}}} \left[T^{\operatorname{curl}}\hat{\psi}_I \circ \mathbf{x}_{\mathcal{K}}(\boldsymbol{\xi})\right]^{\mathsf{T}} \left[T^{\operatorname{curl}}\hat{\psi}_J \circ \mathbf{x}_{\mathcal{K}}(\boldsymbol{\xi})\right] |\mathcal{J}(\boldsymbol{\xi})| d^3\boldsymbol{\xi} + \\
&\quad \int_{\mathcal{K}} \left[\operatorname{curl}\left(T^{\operatorname{curl}}\hat{\psi}_I\right) \circ \mathbf{x}_{\mathcal{K}}(\boldsymbol{\xi})\right]^{\mathsf{T}} \left[\operatorname{curl}\left(T^{\operatorname{curl}}\hat{\psi}_J\right) \circ \mathbf{x}_{\mathcal{K}}(\boldsymbol{\xi})\right] |\mathcal{J}(\boldsymbol{\xi})| d^3\boldsymbol{\xi} \\
&= \int_{\hat{\mathcal{K}}} \left[\left(\mathcal{J}^{-T}\hat{\psi}_I \circ \mathbf{x}_{\mathcal{K}}^{-1}\right) \circ \mathbf{x}_{\mathcal{K}}(\boldsymbol{\xi})\right]^{\mathsf{T}} \left[\left(\mathcal{J}^{-T}\hat{\psi}_J \circ \mathbf{x}_{\mathcal{K}}^{-1}\right) \circ \mathbf{x}_{\mathcal{K}}(\boldsymbol{\xi})\right] |\mathcal{J}(\boldsymbol{\xi})| d^3\boldsymbol{\xi} + \\
&\quad \int_{\mathcal{K}} \left[\left(T^{\operatorname{div}}\widehat{\operatorname{curl}}\hat{\psi}_I\right) \circ \mathbf{x}_{\mathcal{K}}(\boldsymbol{\xi})\right]^{\mathsf{T}} \left[\left(T^{\operatorname{div}}\widehat{\operatorname{curl}}\hat{\psi}_J\right) \circ \mathbf{x}_{\mathcal{K}}(\boldsymbol{\xi})\right] |\mathcal{J}(\boldsymbol{\xi})| d^3\boldsymbol{\xi} \\
&= \int_{\hat{\mathcal{K}}} \hat{\psi}_I(\boldsymbol{\xi})^{\mathsf{T}} \mathcal{D}(\boldsymbol{\xi}) \hat{\psi}_J(\boldsymbol{\xi}) |\mathcal{J}(\boldsymbol{\xi})| d^3\boldsymbol{\xi} + \\
&\quad \int_{\mathcal{K}} \left[\left((|\mathcal{J}|^{-1}\mathcal{J}\widehat{\operatorname{curl}}\hat{\psi}_I) \circ \mathbf{x}_{\mathcal{K}}^{-1}\right) \circ \mathbf{x}_{\mathcal{K}}(\boldsymbol{\xi})\right]^{\mathsf{T}} \left[\left((|\mathcal{J}|^{-1}\mathcal{J}\widehat{\operatorname{curl}}\hat{\psi}_J) \circ \mathbf{x}_{\mathcal{K}}^{-1}\right) \circ \mathbf{x}_{\mathcal{K}}(\boldsymbol{\xi})\right] |\mathcal{J}(\boldsymbol{\xi})| d^3\boldsymbol{\xi} \\
&= \int_{\hat{\mathcal{K}}} \hat{\psi}_I(\boldsymbol{\xi})^{\mathsf{T}} \mathcal{D}(\boldsymbol{\xi}) \hat{\psi}_J(\boldsymbol{\xi}) |\mathcal{J}(\boldsymbol{\xi})| d^3\boldsymbol{\xi} + \int_{\mathcal{K}} \left[\widehat{\operatorname{curl}}\hat{\psi}_I(\boldsymbol{\xi})\right]^{\mathsf{T}} \mathcal{C}(\boldsymbol{\xi}) \left[\widehat{\operatorname{curl}}\hat{\psi}_J(\boldsymbol{\xi})\right] |\mathcal{J}(\boldsymbol{\xi})|^{-1} d^3\boldsymbol{\xi} \\
&= \int_0^1 \int_0^1 \int_0^1 \left\{\hat{\psi}_I(\boldsymbol{\xi})^{\mathsf{T}} \mathcal{D}(\boldsymbol{\xi}) \hat{\psi}_J(\boldsymbol{\xi}) |\mathcal{J}(\boldsymbol{\xi})| + \left[\widehat{\operatorname{curl}}\hat{\psi}_I(\boldsymbol{\xi})\right]^{\mathsf{T}} \mathcal{C}(\boldsymbol{\xi}) \left[\widehat{\operatorname{curl}}\hat{\psi}_J(\boldsymbol{\xi})\right] |\mathcal{J}(\boldsymbol{\xi})|^{-1}\right\} d\xi_3 d\xi_2 d\xi_1 \tag{2.68}
\end{aligned}$$

With the assumption that a 3D finite element code enabled for $H(\text{curl})$ shape functions has the $\widehat{\text{curl}}$ operation incorporated, the conventional algorithm for (2.68) is presented as Algorithm 8.

**Algorithm 8** Conventional computation of the $H(\text{curl})$ Gram Matrix

 **procedure** HCURLGRAM($mid, G^{\text{curl}}$)    ▷ Compute matrix $G^{\text{curl}}$ for element No. $mid$ - Conventional algorithm
  **call** setquadrature3D($mid, p_1, p_2, p_3, L, \boldsymbol{\xi}_{lmn}, w_{lmn}$)
  $G^{\text{curl}} \leftarrow 0$    ▷ Initialize Gram Matrix
  **for** $l, m, n = 1$ to $L$ **do**
   **call** Shape3Hcurl($\boldsymbol{\xi}_{lmn}, p_1, p_2, p_3, \{\hat{\psi}_I(\boldsymbol{\xi}_{lmn})\}, \{\widehat{\text{curl}}\hat{\psi}_I(\boldsymbol{\xi}_{lmn})\}$)  ▷ 3D shape functions at $\boldsymbol{\xi}_{lmn}$
   **call** geometry( $\boldsymbol{\xi}_{lmn}, mid, \mathbf{x}, \mathcal{J}(\boldsymbol{\xi}_{lmn}), \mathcal{J}^{-1}(\boldsymbol{\xi}_{lmn}), |\mathcal{J}|$)  ▷ Compute jacobian
   $\mathcal{D} \leftarrow \mathcal{J}^{-1}(\boldsymbol{\xi}_{lmn})\mathcal{J}^{-T}(\boldsymbol{\xi}_{lmn})$
   $\mathcal{C} \leftarrow \mathcal{J}^{\mathsf{T}}(\boldsymbol{\xi}_{lmn})\mathcal{J}(\boldsymbol{\xi}_{lmn})$
   **for** $J = 0$ to $\dim Q^p - 1$ **do**
    **for** $I = J$ to $\dim Q^p - 1$ **do**
     $G^{\text{curl}}_{IJ} \leftarrow G^{\text{curl}}_{IJ} + \left\{ \hat{\psi}_I(\boldsymbol{\xi}_{lmn})^{\mathsf{T}} \mathcal{D} \hat{\psi}_J(\boldsymbol{\xi}_{lmn})|\mathcal{J}| + \left[\widehat{\text{curl}}\hat{\psi}_I(\boldsymbol{\xi}_{lmn})\right]^{\mathsf{T}} \mathcal{C} \left[\widehat{\text{curl}}\hat{\vartheta}_J(\boldsymbol{\xi}_{lmn})\right] |\mathcal{J}|^{-1} \right\} w_{lmn}$
  **return** $G^{\text{curl}}$

Now we proceed to define a tensor-product shape function for this energy space. As $\hat{\psi}_I, \hat{\psi}_I \in \hat{Q}^p = Q^{p_1-1,p_2,p_3}(\mathcal{I}^3) \times Q^{p_1,p_2-1,p_3}(\mathcal{I}^3) \times Q^{p_1,p_2,p_3-1}(\mathcal{I}^3)$, it follows that:

$$\hat{\psi}_I(\xi_1, \xi_2, \xi_3) = \rho_{a,1;i_1}(\xi_1)\rho_{a,2;i_2}(\xi_2)\rho_{a,3;i_3}(\xi_3)\hat{\mathbf{e}}_a$$
$$\hat{\psi}_J(\xi_1, \xi_2, \xi_3) = \rho_{b,1;j_1}(\xi_1)\rho_{b,2;j_2}(\xi_2)\rho_{b,3;j_3}(\xi_3)\hat{\mathbf{e}}_b. \tag{2.69}$$

where $\hat{\mathbf{e}}_a$, $a = 1, 2, 3$, is the $a$-th canonical cartesian unit vector in the master space; the univariate polynomials $\mu_{a,d;i_d}(\xi_d)$, $d = 1, 2, 3$, are determined through the rule

$$\rho_{a,d;i_d} = \begin{cases} \nu_{i_d} & \text{if } a = d, \\ \chi_{i_d} & \text{otherwise.} \end{cases} \tag{2.70}$$

Also here the values that the integer indices $i_a, j_a$ may take depend on the vector component where $\hat{\psi}_I$ lies ($a$). The correspondence formula between the tensor-product index $I$ and the one-dimensional ones could be of the form

$$I = \begin{cases} i_1 + p_1 i_2 + p_1(p_2+1)i_3 & \text{if a=1,} \\ p_1(p_2+1)(p_3+1) + i_1 + (p_1+1)i_2 + (p_1+1)p_2 i_3 & \text{if a=2,} \\ p_1(p_2+1)(p_3+1) + (p_1+1)p_2(p_3+1) + i_1 + (p_1+1)i_2 + (p_1+1)(p_2+1)i_3 & \text{if a=3.} \end{cases} \tag{2.71}$$

Perhaps it becomes clear by rewriting $\hat{\psi}_I$ using (2.70) just like it was done in the $H(\text{div})$ problem,

$$\hat{\psi}_I(\xi_1, \xi_2, \xi_3) = \nu_{i_a}(\xi_a)\chi_{i_{a+1}}(\xi_{a+1})\chi_{i_{a+2}}(\xi_{a+2})\hat{\mathbf{e}}_a$$
$$= \chi'_{i_{a+1}}(\xi_a)\chi_{i_{a+1}}(\xi_{a+1})\chi_{i_{a+2}}(\xi_{a+2})\hat{\mathbf{e}}_a \tag{2.72}$$

where the relation between $\nu$ and $\chi$ (2.54) was applied. Moreover, the curl operator is the most complicated case of all in this section. Firstly, instead of using explicit indices 1,2,3, we can write the result of this operator having as a reference the component index $a$ in $\hat{\psi}_I$ and making use of the cyclic permutation concept shown above.

$$\widehat{\text{curl}}\hat{\psi}_I = \left[\partial_{a+1}(\hat{\psi}_I)_{a+2} - \partial_{a+2}(\hat{\psi}_I)_{a+1}\right]\hat{\mathbf{e}}_a + \left[\partial_{a+2}(\hat{\psi}_I)_a - \partial_a(\hat{\psi}_I)_{a+2}\right]\hat{\mathbf{e}}_{a+1} + \left[\partial_a(\hat{\psi}_I)_{a+1} - \partial_{a+1}(\hat{\psi}_I)_a\right]\hat{\mathbf{e}}_{a+2}$$
$$= \partial_{a+2}(\hat{\psi}_I)_a\hat{\mathbf{e}}_{a+1} - \partial_{a+1}(\hat{\psi}_I)_a\hat{\mathbf{e}}_{a+2}$$
$$= \chi'_{i_{a+1}}\chi_{i_{a+1}}\chi'_{i_{a+2}}\hat{\mathbf{e}}_{a+1} - \chi'_{i_{a+1}}\chi'_{i_{a+1}}\chi_{i_{a+2}}\hat{\mathbf{e}}_{a+2} \tag{2.73}$$

Now, using the Kronecker delta as a means to manage the superscript (2.37) and the subindex where we need $i_a + 1$ instead of just $i_a$, and introducing this into (2.72,2.73), we can write $\hat{\psi}_I$ and its curl in an appropriate form to perform the nesting of integrals of (2.68).

$$\hat{\psi}_I(\xi_1,\xi_2,\xi_3) = \chi^{\langle\delta_{1a}\rangle}_{i_1+\delta_{1a}}(\xi_1)\chi^{\langle\delta_{2a}\rangle}_{i_2+\delta_{2a}}(\xi_2)\chi^{\langle\delta_{3a}\rangle}_{i_3+\delta_{3a}}(\xi_3)\hat{\mathbf{e}}_a, \tag{2.74}$$

$$\widehat{\text{curl}}\hat{\psi}_I(\xi_1,\xi_2,\xi_3) = \chi^{\langle 1-\delta_{1(a+1)}\rangle}_{i_1+\delta_{1a}}(\xi_1)\chi^{\langle 1-\delta_{2(a+1)}\rangle}_{i_2+\delta_{2a}}(\xi_2)\chi^{\langle 1-\delta_{3(a+1)}\rangle}_{i_3+\delta_{3a}}(\xi_3)\hat{\mathbf{e}}_{a+1} -$$
$$\chi^{\langle 1-\delta_{1(a+2)}\rangle}_{i_1+\delta_{1a}}(\xi_1)\chi^{\langle 1-\delta_{2(a+2)}\rangle}_{i_2+\delta_{2a}}(\xi_2)\chi^{\langle 1-\delta_{3(a+2)}\rangle}_{i_3+\delta_{3a}}(\xi_3)\hat{\mathbf{e}}_{a+2}. \tag{2.75}$$

By applying (2.74) and (2.75) to the first and second terms of (2.68), respectively, we converge to a definite expression that we can finally implement for the computation of the Gram matrix entry $G^{\text{curl}}_{IJ}$. Putting aside momentarily the $|\mathcal{J}|$ factors we have the following equations for each term,

$$\hat{\psi}_I(\boldsymbol{\xi})^\mathsf{T}\mathcal{D}(\boldsymbol{\xi})\hat{\psi}_J(\boldsymbol{\xi}) = \chi^{\langle\delta_{1a}\rangle}_{i_1+\delta_{1a}}(\xi_1)\chi^{\langle\delta_{1b}\rangle}_{j_1+\delta_{1b}}(\xi_1)\chi^{\langle\delta_{2a}\rangle}_{i_2+\delta_{2a}}(\xi_2)\chi^{\langle\delta_{2b}\rangle}_{j_2+\delta_{2b}}(\xi_2)\chi^{\langle\delta_{3a}\rangle}_{i_3+\delta_{3a}}(\xi_3)\chi^{\langle\delta_{3b}\rangle}_{j_3+\delta_{3b}}(\xi_3)\hat{\mathbf{e}}_a^\mathsf{T}\mathcal{D}(\boldsymbol{\xi})\hat{\mathbf{e}}_b, \tag{2.76}$$

$$\left[\widehat{\text{curl}}\hat{\psi}_I\right]^\mathsf{T}\mathcal{C}\left[\widehat{\text{curl}}\hat{\psi}_J\right] = \chi^{\langle 1-\delta_{1(a+1)}\rangle}_{i_1+\delta_{1a}}\chi^{\langle 1-\delta_{1(b+1)}\rangle}_{j_1+\delta_{1b}}\chi^{\langle 1-\delta_{2(a+1)}\rangle}_{i_2+\delta_{2a}}\chi^{\langle 1-\delta_{2(b+1)}\rangle}_{j_2+\delta_{2b}}\chi^{\langle 1-\delta_{3(a+1)}\rangle}_{i_3+\delta_{3a}}\chi^{\langle 1-\delta_{3(b+1)}\rangle}_{j_3+\delta_{3b}}\hat{\mathbf{e}}_{a+1}^\mathsf{T}\mathcal{C}\hat{\mathbf{e}}_{b+1} -$$
$$\chi^{\langle 1-\delta_{1(a+2)}\rangle}_{i_1+\delta_{1a}}\chi^{\langle 1-\delta_{1(b+1)}\rangle}_{j_1+\delta_{1b}}\chi^{\langle 1-\delta_{2(a+2)}\rangle}_{i_2+\delta_{2a}}\chi^{\langle 1-\delta_{2(b+1)}\rangle}_{j_2+\delta_{2b}}\chi^{\langle 1-\delta_{3(a+2)}\rangle}_{i_3+\delta_{3a}}\chi^{\langle 1-\delta_{3(b+1)}\rangle}_{j_3+\delta_{3b}}\hat{\mathbf{e}}_{a+2}^\mathsf{T}\mathcal{C}\hat{\mathbf{e}}_{b+1} -$$
$$\chi^{\langle 1-\delta_{1(a+1)}\rangle}_{i_1+\delta_{1a}}\chi^{\langle 1-\delta_{1(b+2)}\rangle}_{j_1+\delta_{1b}}\chi^{\langle 1-\delta_{2(a+1)}\rangle}_{i_2+\delta_{2a}}\chi^{\langle 1-\delta_{2(b+2)}\rangle}_{j_2+\delta_{2b}}\chi^{\langle 1-\delta_{3(a+1)}\rangle}_{i_3+\delta_{3a}}\chi^{\langle 1-\delta_{3(b+2)}\rangle}_{j_3+\delta_{3b}}\hat{\mathbf{e}}_{a+1}^\mathsf{T}\mathcal{C}\hat{\mathbf{e}}_{b+2} +$$
$$\chi^{\langle 1-\delta_{1(a+2)}\rangle}_{i_1+\delta_{1a}}\chi^{\langle 1-\delta_{1(b+2)}\rangle}_{j_1+\delta_{1b}}\chi^{\langle 1-\delta_{2(a+2)}\rangle}_{i_2+\delta_{2a}}\chi^{\langle 1-\delta_{2(b+2)}\rangle}_{j_2+\delta_{2b}}\chi^{\langle 1-\delta_{3(a+2)}\rangle}_{i_3+\delta_{3a}}\chi^{\langle 1-\delta_{3(b+2)}\rangle}_{j_3+\delta_{3b}}\hat{\mathbf{e}}_{a+2}^\mathsf{T}\mathcal{C}\hat{\mathbf{e}}_{b+2}. \tag{2.77}$$

Note that in the second equation the arguments were omitted for practical reasons, but it is easy to determine what is the argument of each univariate polynomial by comparing to (2.76). We can further recall that $\hat{\mathbf{e}}_a^\mathsf{T}\mathcal{D}(\boldsymbol{\xi})\hat{\mathbf{e}}_b = D_{ab}$, $\hat{\mathbf{e}}_{a+1}^\mathsf{T}\mathcal{C}\hat{\mathbf{e}}_{b+1} = C_{(a+1)(b+1)}$, and similarly we can retrieve $C_{(a+2)(b+1)}, C_{(a+1)(b+2)}, C_{(a+2)(b+2)}$. The second term of the Gram matrix is more intricate as it comprises four subterms, each of which needs to be computed separately. However, it is possible to generalize the definition of the auxiliary functions adding a new pair of indices. Then, let

$$\mathcal{G}^{\text{curl } A;\alpha\beta}_{ab;i_3j_3}(\xi_1,\xi_2) = \int_0^1 \chi^{\langle 1-\delta_{3(a+\alpha)}\rangle}_{i_3+\delta_{3a}}(\xi_3)\chi^{\langle 1-\delta_{3(b+\beta)}\rangle}_{j_3+\delta_{3b}}(\xi_3)C_{(a+\alpha)(b+\beta)}(\xi_1,\xi_2,\xi_3)|\mathcal{J}(\xi_1,\xi_2,\xi_3)|^{-1}d\xi_3, \tag{2.78}$$

$$\mathcal{G}^{\text{curl } B;\alpha\beta}_{ab;i_2j_2i_3j_3}(\xi_1) = \int_0^1 \chi^{\langle 1-\delta_{2(a+\alpha)}\rangle}_{i_2+\delta_{2a}}(\xi_2)\chi^{\langle 1-\delta_{2(b+\beta)}\rangle}_{j_2+\delta_{2b}}(\xi_2)\mathcal{G}^{\text{curl } A;\alpha\beta}_{ab;i_3j_3}(\xi_1,\xi_2)d\xi_2, \tag{2.79}$$

$$\mathcal{G}^{\text{curl};\alpha\beta}_{ab;i_1j_1i_2j_2i_3j_3} = \int_0^1 \chi^{\langle 1-\delta_{1(a+\alpha)}\rangle}_{i_1+\delta_{1a}}(\xi_1)\chi^{\langle 1-\delta_{1(b+\beta)}\rangle}_{j_1+\delta_{1b}}(\xi_1)\mathcal{G}^{\text{curl } B;\alpha\beta}_{ab;i_2j_2i_3j_3}(\xi_1)d\xi_1, \tag{2.80}$$

where $\alpha, \beta = 1, 2$. In the same fashion, for the first term the auxiliary functions are

$$\check{\mathcal{G}}^{\text{curl }A}_{ab;i_3j_3}(\xi_1, \xi_2) = \int_0^1 \chi^{\langle\delta_{3a}\rangle}_{i_3+\delta_{3a}}(\xi_3)\chi^{\langle\delta_{3b}\rangle}_{j_3+\delta_{3b}}(\xi_3)D_{ab}(\xi_1,\xi_2,\xi_3)|\mathcal{J}(\xi_1,\xi_2,\xi_3)|d\xi_3, \quad (2.81)$$

$$\check{\mathcal{G}}^{\text{curl }B}_{ab;i_2j_2i_3j_3}(\xi_1) = \int_0^1 \chi^{\langle\delta_{2a}\rangle}_{i_2+\delta_{2a}}(\xi_2)\chi^{\langle\delta_{2b}\rangle}_{j_2+\delta_{2b}}(\xi_2)\check{\mathcal{G}}^{\text{curl }A}_{ab;i_3j_3}(\xi_1,\xi_2)d\xi_2, \quad (2.82)$$

$$\check{\mathcal{G}}^{\text{curl}}_{ab;i_1j_1i_2j_2i_3j_3} = \int_0^1 \chi^{\langle\delta_{1a}\rangle}_{i_1+\delta_{1a}}(\xi_1)\chi^{\langle\delta_{1a}\rangle}_{j_1+\delta_{1b}}(\xi_1)\check{\mathcal{G}}^{\text{curl }B}_{ab;i_2j_2i_3j_3}(\xi_1)d\xi_1, \quad (2.83)$$

Finally, every Gram matrix entry is computed as

$$G^{\text{curl}}_{IJ} = \check{\mathcal{G}}^{\text{curl}}_{ab;i_1j_1i_2j_2i_3j_3} + \sum_{\alpha,\beta=1}^{2}(-1)^{\alpha+\beta}\mathcal{G}^{\text{curl};\alpha\beta}_{ab;i_1j_1i_2j_2i_3j_3}. \quad (2.84)$$

Algorithm 9 implements (2.78)-(2.84) using one-dimensional numerical integration for each auxiliary function. The symmetries are not explicitly given, but if the reader wishes to investigate them and incorporate them, those could be applied within a function that here is called $fillsymmetric3Hcurl$, which would additionally perform the sum of the two terms of this inner product (implementing this kind of function is not required for the algorithm to work). If we have a uniform-degree polynomial space (in the sense of the exact sequence), then this algorithm has a complexity of $\mathcal{O}[p^3(p+1)^4] \sim \mathcal{O}(p^7)$ whereas Algorithm 8 accounts for $\mathcal{O}[p^5(p+1)^4] \sim \mathcal{O}(p^9)$, assuming we use $L=p$ quadrature points.

## 2.7 Vector- and matrix-valued versions

In the former subsection all the algorithms for the different Gram matrices have been introduced. However, some variational formulations require higher-dimension versions of the basic spaces so far presented. For instance, the two scalar-valued function spaces treated above, $L^2(\mathcal{K})$ and $H^1(\mathcal{K})$, have their vector-valued versions $\boldsymbol{L}^2(\mathcal{K}) := (L^2(\mathcal{K}))^3$ and $\boldsymbol{H}^1(\mathcal{K}) := (H^1(\mathcal{K}))^3$, while the two vector-valued energy spaces $H(\text{div},\mathcal{K})$ and $H(\text{curl},\mathcal{K})$ have matrix-valued counterparts, $\boldsymbol{H}(\text{div},\mathcal{K}) := (H(\text{div},\mathcal{K}))^3$ and $\boldsymbol{H}(\text{curl},\mathcal{K}) := (H(\text{curl},\mathcal{K}))^3$ which essentially consist of the space of matrices with each row being a vector in the original space. As a sample to each case (scalar- and vector-valued) let us present $\boldsymbol{H}^1(\mathcal{K})$ and $\boldsymbol{H}(\text{div},\mathcal{K})$.

Let $\Phi_I, \Phi_J \in \boldsymbol{H}^1(\mathcal{K})$ such that

$$\Phi_I = \varphi_i \mathbf{e}_a,$$
$$\Phi_J = \varphi_j \mathbf{e}_b,$$

where $a, b = 1, 2, 3$; $\mathbf{e}_a, \mathbf{e}_b$ are canonical cartesian unit vectors in the physical space, and $\varphi_i, \varphi_j \in H^1(\mathcal{K})$.

The inner product between those functions is

$$(\Phi_I, \Phi_J)_{\boldsymbol{H}^1(\mathcal{K})} := \int_\mathcal{K} \Phi_I^\mathsf{T}\Phi_J d^3\boldsymbol{x} + \int_\mathcal{K} (\nabla\Phi_I):(\nabla\Phi_J)\,d^3\boldsymbol{x}, \quad (2.85)$$

$$= \int_\mathcal{K} \varphi_i \mathbf{e}_a^\mathsf{T}\mathbf{e}_b\varphi_j d^3\boldsymbol{x} + \int_\mathcal{K} (\mathbf{e}_a \otimes \nabla\varphi_i):(\mathbf{e}_b \otimes \nabla\varphi_j)\,d^3\boldsymbol{x},$$

**Algorithm 9** Tensorization-based computation of the $H(\mathrm{curl})$ Gram Matrix
---
  **procedure** HCURLGRAMTENSOR($mid, G^{\mathrm{curl}}$)          ▷ $G^{\mathrm{curl}}$ for element No. $mid$ -Tensorization - based algorithm
    $p_{\max} \leftarrow \max\{p_1, p_2, p_3\}$
    **call** setquadrature1D($mid, p_{\max}, L, \zeta^l, w^l$)
    $\check{\mathcal{G}}^{\mathrm{grad}}, \mathcal{G}^{\mathrm{grad}} \leftarrow 0$          ▷ Initialize Gram Matrix
    **for** $l = 1$ to $L$ **do**
        **call** Shape1H1($\zeta^l, p_1, \chi_{i_1}(\zeta^l), \chi'_{i_1}(\zeta^l)\}$)         ▷ Evaluate 1D shape functions at $\zeta^l$
        **for** $j_3, i_3 = 0$ to $p_3$ **do**
            $\check{\mathcal{G}}^{\mathrm{curl}\,B}, \mathcal{G}^{\mathrm{curl}\,B} \leftarrow 0$
            **for** $m = 1$ to $L$ **do**
                **call** Shape1H1($\zeta^m, p_2, \chi_{i_2}(\zeta^m), \chi'_{i_2}(\zeta^m)\}$)     ▷ Evaluate 1D shape functions at $\zeta^m$
                $\check{\mathcal{G}}^{\mathrm{curl}\,A}, \mathcal{G}^{\mathrm{curl}\,A} \leftarrow 0$
                **for** $n = 1$ to $L$ **do**
                    **call** Shape1H1($\zeta^n, p_3, \chi_{i_3}(\zeta^n), \chi'_{i_3}(\zeta^n)\}$)     ▷ Evaluate 1D shape functions at $\zeta^n$
                    $\boldsymbol{\xi}_{lmn} \leftarrow (\zeta^l, \zeta^m, \zeta^n)$
                    **call** geometry($\boldsymbol{\xi}_{lmn}, mid, \mathbf{x}, \mathcal{J}(\boldsymbol{\xi}_{lmn}), \mathcal{J}^{-1}(\boldsymbol{\xi}_{lmn}), |\mathcal{J}|$)   ▷ Compute jacobian
                    $\mathcal{D} \leftarrow \mathcal{J}^{-1}(\boldsymbol{\xi}_{lmn})\mathcal{J}^{-T}(\boldsymbol{\xi}_{lmn})$
                    $\mathcal{C} \leftarrow \mathcal{J}^{\mathsf{T}}(\boldsymbol{\xi}_{lmn})\mathcal{J}(\boldsymbol{\xi}_{lmn})$
                    **for** $a, b = 1$ to $3$ **do**
                        **if** $j_3 + \delta_{3b} \leq p_3$ and $i_3 + \delta_{3a} \leq p_3$ **then**         ▷ Avoids extra computations
                            $\check{\mathcal{G}}^A_{ab;i_3j_3} \leftarrow \check{\mathcal{G}}^A_{ab;i_3j_3} + \chi^{\langle\delta_{3a}\rangle}_{i_3+\delta_{3a}}(\zeta^n)\chi^{\langle\delta_{3b}\rangle}_{j_3+\delta_{3b}}(\zeta^n)D_{ab}|\mathcal{J}|w^n$
                            **for** $\alpha, \beta = 1$ to $2$ **do**
                                $\mathcal{G}^{\mathrm{curl}\,A;\alpha\beta}_{ab;i_3j_3} \leftarrow \mathcal{G}^{\mathrm{curl}\,A;\alpha\beta}_{ab;i_3j_3} + \chi^{\langle 1-\delta_{3a}\rangle}_{i_3+\delta_{3a}}(\zeta^n)\chi^{\langle 1-\delta_{3b}\rangle}_{j_3+\delta_{3b}}(\zeta^n)C_{(a+\alpha)(b+\beta)}|\mathcal{J}|^{-1}w^n$   ▷ (2.78)
            **for** $j_2, i_2 = 0$ to $p_2$ **do**
                **for** $a, b = 1$ to $3$ **do**
                    **if** $j_2 + \delta_{2b} \leq p_2$ and $i_2 + \delta_{2a} \leq p_2$ **then**         ▷ Avoids extra computations
                        $\check{\mathcal{G}}^B_{ab;i_2j_2i_3j_3} \leftarrow \check{\mathcal{G}}^B_{ab;i_2j_2i_3j_3} + \chi^{\langle\delta_{2a}\rangle}_{i_2+\delta_{2a}}(\zeta^m)\chi^{\langle\delta_{2b}\rangle}_{j_2+\delta_{2b}}(\zeta^m)\check{\mathcal{G}}^A_{ab;i_3j_3}w^m$
                        **for** $\alpha, \beta = 1$ to $2$ **do**
                            $\mathcal{G}^{\mathrm{curl}\,B;\alpha\beta}_{ab;i_2j_2i_3j_3} \leftarrow \mathcal{G}^{\mathrm{curl}\,B;\alpha\beta}_{ab;i_2j_2i_3j_3} + \chi^{\langle 1-\delta_{2a}\rangle}_{i_2+\delta_{2a}}(\zeta^m)\chi^{\langle 1-\delta_{2b}\rangle}_{j_2+\delta_{2b}}(\zeta^m)\mathcal{G}^{\mathrm{curl}\,A;\alpha\beta}_{ab;i_3j_3}w^m$   ▷ (2.79)
        **for** $j_2, i_2 = 0$ to $p_2$ **do**
            **for** $j_1, i_1 = 0$ to $p_1$ **do**
                Determine $I, J$ through (2.71)
                **if** $J \geq I$ **then**
                    **for** $a, b = 1$ to $3$ **do**
                        **if** $j_1 + \delta_{1b} \leq p_1$ and $i_1 + \delta_{1a} \leq p_1$ **then**         ▷ Avoids extra computations
                            $\check{\mathcal{G}}_{ab;i_1j_1i_2j_2i_3j_3} \leftarrow \check{\mathcal{G}}_{ab;i_1j_1i_2j_2i_3j_3} + \chi^{\langle\delta_{1a}\rangle}_{i_1+\delta_{1a}}(\zeta^l)\chi^{\langle\delta_{1b}\rangle}_{j_1+\delta_{1b}}(\zeta^l)\check{\mathcal{G}}^B_{ab;i_2j_2i_3j_3}w^l$
                            **for** $\alpha, \beta = 1$ to $2$ **do**
                                $\mathcal{G}^{\mathrm{curl};\alpha\beta}_{ab;i_1j_1i_2j_2i_3j_3} \leftarrow \mathcal{G}^{\mathrm{curl};\alpha\beta}_{ab;i_1j_1i_2j_2i_3j_3} + \chi^{\langle 1-\delta_{1a}\rangle}_{i_1+\delta_{1a}}(\zeta^l)\chi^{\langle 1-\delta_{1b}\rangle}_{j_1+\delta_{1b}}(\zeta^l)\mathcal{G}^{\mathrm{curl}\,B;\alpha\beta}_{ab;i_2j_2i_3j_3}w^l$
    **call** fillsymmetric3Hcurl($p_1, p_2, p_3, \check{\mathcal{G}}, \mathcal{G}^{\mathrm{curl}}, G^{\mathrm{curl}}$)         ▷ Adds terms and fills matrix using symmetries
    **return** $G^{\mathrm{curl}}$

$$=\delta_{ab}\int_{\mathcal{K}}\varphi_i\varphi_j d^3\boldsymbol{x}+\delta_{ab}\int_{\mathcal{K}}[\nabla\varphi_i]^{\mathsf{T}}\nabla\varphi_j d^3\boldsymbol{x},$$
$$=\delta_{ab}\left(\varphi_i,\varphi_j\right)_{H^1(\mathcal{K})}. \tag{2.86}$$

Similarly, for $\Theta_I, \Theta_J \in \boldsymbol{H}(\mathrm{div},\mathcal{K})$, such that

$$\Theta_I = \mathbf{e}_a \otimes \vartheta_i,$$
$$\Theta_J = \mathbf{e}_b \otimes \vartheta_j,$$

where $a, b = 1, 2, 3$, and $\vartheta_i, \vartheta_j \in H(\mathrm{div},\mathcal{K})$. Their inner product yields

$$(\Theta_I,\Theta_J)_{\boldsymbol{H}(\mathrm{div},\mathcal{K})} := \int_{\mathcal{K}} \Theta_I : \Theta_J d^3\boldsymbol{x} + \int_{\mathcal{K}} (\mathrm{div}\,\Theta_I)^{\mathsf{T}}\,(\mathrm{div}\,\Theta_J)\, d^3\boldsymbol{x}, \tag{2.87}$$
$$= \int_{\mathcal{K}} (\mathbf{e}_a \otimes \vartheta_i) : (\mathbf{e}_b \otimes \vartheta_j)\, d^3\boldsymbol{x} + \int_{\mathcal{K}} (\mathrm{div}\,\vartheta_i \mathbf{e}_a)^{\mathsf{T}}\,(\mathrm{div}\,\vartheta_j \mathbf{e}_b)\, d^3\boldsymbol{x},$$
$$=\delta_{ab}\int_{\mathcal{K}}\vartheta_i\vartheta_j d^3\boldsymbol{x}+\delta_{ab}\int_{\mathcal{K}}\mathrm{div}\,\vartheta_i\,\mathrm{div}\,\vartheta_j d^3\boldsymbol{x},$$
$$=\delta_{ab}\left(\vartheta_i,\vartheta_j\right)_{H(\mathrm{div},\mathcal{K})}. \tag{2.88}$$

Notice that (2.86) and (2.88) suggest that, in order to compute a Gram matrix for these higher-dimension spaces, the calculation of each entry is equal to the associated Gram matrix entry in the original space, reproduced three times, i.e. when $a = b$. The indices $a, b$ may be introduced as new loops added within the algorithms or we could opt to keep these identical and just change the filling functions so that the new Gram matrix replicates three times the original one. This consideration should be implemented when the problems like those in the Applications section require this kind of test spaces. The same result is extensive to $\boldsymbol{L}^2(\mathcal{K})$ and $\boldsymbol{H}(\mathrm{curl},\mathcal{K})$.

## 2.8 Practical simplifications

In this part of the article we want to focus on those very frequent cases where the element map is an affine map or when the physical element coincides with a protrusion of a plane 2D shape, and the coordinate associated to the depth (e.g. $x_1$) is given by an affine function of only one of the parameters, say $\xi_1$. In other words, if we have a Jacobian coinciding with either of the forms

$$\mathcal{J} = \begin{pmatrix} \lambda_1 & 0 & 0 \\ 0 & \lambda_2 & 0 \\ 0 & 0 & \lambda_3 \end{pmatrix} \quad \text{or} \quad \mathcal{J}(\xi_2,\xi_3) = \begin{pmatrix} \lambda_1 & 0 & 0 \\ 0 & \eta_{22}(\xi_2,\xi_3) & \eta_{23}(\xi_2,\xi_3) \\ 0 & \eta_{32}(\xi_2,\xi_3) & \eta_{33}(\xi_2,\xi_3) \end{pmatrix}, \tag{2.89}$$

where $\lambda_1, \lambda_2, \lambda_3$ are nonzero constants and $\eta_{ab} := \partial \mathbf{x}_a / \partial \xi_b$, then we can exploit such a fact and simplify our algorithms. In either situation of (2.89) the algorithms may be led to one order of magnitude cheaper by removing one or more loops in their structures and picking very appropriate shape functions for which we can load simple formulas of their integrals.

### 2.8.1 Simpler parametric mappings

In the first case shown in (2.89) we have as a result that both $|\mathcal{J}|$ and $\mathcal{D}$ are constant (therefore $\mathcal{C} = \mathcal{D}^{-1}$ too). It follows that we can pull them out of the integrals and in the end (taking $L^2$ as an example may illustrate the idea for the remaining spaces), we have any Gram matrix entry (2.19) transformed into a much simpler expression,

$$\begin{aligned}(v_I, v_J) &= \int_0^1 \nu_{i_1}(\xi_1)\nu_{j_1}(\xi_1) \left\{ \int_0^1 \nu_{i_2}(\xi_2)\nu_{j_2}(\xi_2) \left[ \int_0^1 \nu_{i_3}(\xi_3)\nu_{j_3}(\xi_3)|\mathcal{J}|^{-1}d\xi_3 \right] d\xi_2 \right\} d\xi_1 \\ &= |\mathcal{J}|^{-1} \int_0^1 \nu_{i_1}(\xi_1)\nu_{j_1}(\xi_1) \left\{ \int_0^1 \nu_{i_2}(\xi_2)\nu_{j_2}(\xi_2) \left[ \int_0^1 \nu_{i_3}(\xi_3)\nu_{j_3}(\xi_3)d\xi_3 \right] d\xi_2 \right\} d\xi_1 \\ &= |\mathcal{J}|^{-1} \int_0^1 \nu_{i_1}(\xi_1)\nu_{j_1}(\xi_1)d\xi_1 \int_0^1 \nu_{i_2}(\xi_2)\nu_{j_2}(\xi_2)d\xi_2 \int_0^1 \nu_{i_3}(\xi_3)\nu_{j_3}(\xi_3)d\xi_3.\end{aligned}$$

Using the relation (2.54) we see that all the equation depends solely on what the basis $\{\chi_i\}$ contains. If we define a new auxiliary function that represents that interval integral of two shape functions (or of their derivatives) then we can achieve a simplified version of the derived formulas. For $c = 1, 2, 3$, $0 \leq i_c, j_c \leq p_c$ and $r, s = 0, 1$ let

$$F^{rs}_{c;i_c j_c} := \int_0^1 \chi^{\langle r \rangle}_{i_c}(\xi_c) \chi^{\langle s \rangle}_{j_c}(\xi_c) d\xi_c. \tag{2.90}$$

If we know the shape function set, then we can have precomputed every $F^{rs}_{c;i_c j_c}$. By introducing this approach into Algorithm 2 we abolish three loops over quadrature points. Algorithm 10 applies this result.

On the other hand, if we have the second situation in (2.89) we still need to integrate over $\xi_2$ and $\xi_3$, but just by keeping the order of the nesting of integrals, we still can reduce the cost by one order of magnitude. Algorithm 11 makes use of this idea and the reader can verify that the ultimate accumulation statement is performed just $\mathcal{O}(p^6)$ times.

The extension of the ideas implemented in Algorithms 10 and 11 are directly extensive to the rest of the spaces studied throughout this section, and for high values of $p$ it really can represent a significant amount of computational savings when going from $\mathcal{O}(p^7)$ to $\mathcal{O}(p^6)$. Nonetheless, this is only possible if some additional effort in coding is carried out, both in the logical structure reorderings we have appreciated but also in having precomputed results for $F^{rs}_{c;i_c j_c}$, which is not a conventional task, especially because FE codes tend to be thought for most generalized problems. Then, it turns pretty helpful to pick shape functions that ease that precomputation load, or that can deliver simple formulas for any case of $F^{rs}_{c;i_c j_c}$. Legendre polynomials satisfy that requirement as in the best of cases (purely constant jacobian) we could even circumvent every shape function evaluation. If we are in the less simplified case, it still adapts to our requirement of being able to compute the integrals of at least one coordinate with no quadrature points. And if there is no possibility of simplification, then Legendre polynomials will at least comply with the hierarchical basis property, making the algorithms developed above directly applicable.

**Algorithm 10** Tensorization-based computation of the $L^2$ Gram Matrix-Simplified 1

**procedure** L2GRAMTENSORS1($mid, G$) ▷ Compute matrix $G$ for element No. $mid$ - Tensorization-based algorithm
    $\mathcal{G} \leftarrow 0$ ▷ Initialize Gram Matrix
    **for** $j_3 = 0$ to $p_3 - 1$ **do**
        **for** $i_3 = j_3$ to $p_3 - 1$ **do**
            $\mathcal{G}^B \leftarrow 0$
            $\mathcal{G}^A \leftarrow 0$
            $\boldsymbol{\xi}_{lmn} \leftarrow (0.5, 0.5, 0.5)$ ▷ Evaluate element map at the midpoint
            **call** geometry($\boldsymbol{\xi}_{lmn}, mid, \mathbf{x}, \mathcal{J}(\boldsymbol{\xi}_{lmn}), \mathcal{J}^{-1}(\boldsymbol{\xi}_{lmn}), |\mathcal{J}|$) ▷ Compute jacobian
            $\mathcal{G}^A_{i_3 j_3} \leftarrow F^{11}_{3;(i_3+1)(j_3+1)} |\mathcal{J}|^{-1}$
            **for** $j_2 = 0$ to $p_2 - 1$ **do**
                **for** $i_2 = j_2$ to $p_2 - 1$ **do**
                    $\mathcal{G}^B_{i_2 j_2 i_3 j_3} \leftarrow F^{11}_{1;(i_2+1)(j_2+1)} \mathcal{G}^A_{i_3 j_3}$
            **for** $j_2 = 0$ to $p_2 - 1$ **do**
                **for** $i_2 = j_2$ to $p_2 - 1$ **do**
                  **for** $j_1 = 0$ to $p_1 - 1$ **do**
                      **for** $i_1 = j_1$ to $p_1 - 1$ **do**
                          $\mathcal{G}_{i_1 j_1 i_2 j_2 i_3 j_3} \leftarrow F^{11}_{1;(i_1+1)(j_1+1)} \mathcal{G}^B_{i_2 j_2 i_3 j_3}$
    **call** fillsymmetric3L2($p_1, p_2, p_3, \mathcal{G}, G$) ▷ Fills matrix using symmetries
    **return** $G$

### 2.8.2 Use of Legendre polynomials

Legendre polynomials are a family of hierarchical and $L^2$ orthogonal polynomials defined over the interval [-1,1], with the average zero property (except the constant corresponding to $0^{\text{th}}$ degree). By shifting the argument we can keep all those properties over the interval [0,1]. We are denoting by $P_i$ the $i^{\text{th}}$-degree Legendre polynomial, with $i = 0, 1, 2, ...$. Just by defining the first two members of the set along with a recursion formula, we can obtain all of the Legendre polynomial family over the closed master interval.

$$P_0(\xi) = 1 \tag{2.91}$$

$$P_1(\xi) = 2\xi - 1 \tag{2.92}$$

$$i P_i(\xi) = (2i-1)(2\xi - 1) P_{i-1}(\xi) - (i-1) P_{i-2}(\xi), \quad i \geq 2. \tag{2.93}$$

In (2.93) the recursion formula is built such that $P_i(1) = 1$ for every non-negative $i$.

Additionally, the mentioned properties of the Legendre polynomials are formally rewritten like this:

1. Hierarchical polynomial basis

$$\mathcal{P}^n([0,1]) = \text{span}\{P_0, ..., P_n\}, \quad n = 0, 1, 2, 3, ... \tag{2.94}$$

**Algorithm 11** Tensorization-based computation of the $L^2$ Gram Matrix-Simplified element map
---
**procedure** L2GRAMTENSORS2($mid, G$) ▷ Compute matrix $G$ for element No. $mid$ - Tensorization-based algorithm
    $p_{\max} \leftarrow \max\{p_1, p_2, p_3\}$
    **call** setquadrature1D($mid, p_{\max} - 1, L, \zeta^l, w^l$)
    $\mathcal{G} \leftarrow 0$     ▷ Initialize Gram Matrix
    **for** $j_3 = 0$ to $p_3 - 1$ **do**
        **for** $i_3 = j_3$ to $p_3 - 1$ **do**
            $\mathcal{G}^A \leftarrow 0$
            **for** $m = 1, L$ **do**
                **call** Shape1L2($\zeta^m, p_2, \nu_{i_2}(\zeta^m)\}$)     ▷ Evaluate 1D shape functions at $\zeta^m$
                $\mathcal{G}^B \leftarrow 0$
                **for** $n = 1$ to $L$ **do**
                      **call** Shape1L2($\zeta^n, p_3, \nu_{i_3}(\zeta^n)\}$)     ▷ Evaluate 1D shape functions at $\zeta^n$
                      $\boldsymbol{\xi}_{lmn} \leftarrow (0.5, \zeta^m, \zeta^n)$     ▷ Use 0.5 instead of any $\zeta^l$
                      **call** geometry($\boldsymbol{\xi}_{lmn}, mid, \mathbf{x}, \mathcal{J}(\boldsymbol{\xi}_{lmn}), \mathcal{J}^{-1}(\boldsymbol{\xi}_{lmn}), |\mathcal{J}|$)     ▷ Compute jacobian
                      $\mathcal{G}^A_{i_3 j_3} \leftarrow \mathcal{G}^A_{i_3 j_3} + \nu_{i_3}(\zeta^n)\nu_{j_3}(\zeta^n)|\mathcal{J}|^{-1}w^n$     ▷ Accumulate through (2.23)
                **for** $j_2 = 0$ to $p_2 - 1$ **do**
                    **for** $i_2 = j_2$ to $p_2 - 1$ **do**
                      $\mathcal{G}^B_{i_2 j_2 i_3 j_3} \leftarrow \mathcal{G}^B_{i_2 j_2 i_3 j_3} + \nu_{i_2}(\zeta^m)\nu_{j_2}(\zeta^m)\mathcal{G}^A_{i_3 j_3}(\zeta^l, \zeta^m)w^m$     ▷ From (2.24)
            **for** $j_2 = 0$ to $p_2 - 1$ **do**
                **for** $i_2 = j_2$ to $p_2 - 1$ **do**
                  **for** $j_1 = 0$ to $p_1 - 1$ **do**
                      **for** $i_1 = j_1$ to $p_1 - 1$ **do**
                          $\mathcal{G}_{i_1 j_1 i_2 j_2 i_3 j_3} \leftarrow F^{11}_{1;(i_1+1)(j_1+1)} \mathcal{G}^B_{i_2 j_2 i_3 j_3}$
    **call** fillsymmetric3L2($p_1, p_2, p_3, \mathcal{G}, G$)     ▷ Fills matrix using symmetries
    **return** $G$

---

2. Orthogonality

$$(P_i, P_j)_{L^2(\mathcal{I})} = \int_0^1 P_i(\xi)P_j(\xi)d\xi = \delta_{ij}\frac{1}{2i+1}, \quad i, j \geq 0, \tag{2.95}$$

with $\delta_{ij}$ being the Kronecker delta.

3. Average zero

$$\int_0^1 P_i(\xi)d\xi = 0, \quad i \geq 1. \tag{2.96}$$

Another useful property derived from (2.91)-(2.93) is that Legendre polynomials of odd degree are odd with respect to the line $\xi = 1/2$, and similarly those of even index are even with respect to that line.

However, we usually need to evaluate derivatives of the shape functions. In that sense, it would be preferable to have all the nice properties of Legendre polynomials in the derivatives of our shape functions. This motivates the introduction of the integrated Legendre polynomials, defined as:

$$L_i(\xi) = \int_0^\xi P_{i-1}(t)dt = 0, \quad i \geq 1 \tag{2.97}$$

This implies:

$$L_i'(\xi) = P_{i-1}(\xi), \quad i \geq 1. \tag{2.98}$$

Notice that (2.96) and (2.97) imply that $L_i(0) = L_i(1) = 0$ for all $i \geq 2$. This means that all higher order integrated Legendre polynomials are bubbles (smooth functions vanishing at the boundary of its support). Now, with some elementary calculus and the definitions and properties above we can find a corresponding recursion formula for this new family of polynomials:

$$L_1(\xi) = \xi \tag{2.99}$$
$$2(2i-1)L_i(\xi) = P_i(\xi) - P_{i-2}(\xi), \quad i \geq 2. \tag{2.100}$$

This recursive definition means that, for $i \geq 2$, the property of being either even or odd with respect to $\xi = 1/2$ correspondingly to their indices still holds. Finally, we can also do a bit of notation abuse and add a new function to this family:

$$L_0(\xi) = 1 - \xi = 1 - L_1(\xi). \tag{2.101}$$

The definition of this additional function has the purpose of completing a full basis with the $L_i$ polynomials that is also hierarchical (as long as $i \geq 1$), just like in the original Legendre family. Thus we have

$$\mathcal{P}^n([0,1]) = \text{span}\{L_0, ..., L_n\}, \quad n = 1, 2, 3, ... \tag{2.102}$$

If the 1D $H^1$ shape functions are integrated Legendre polynomials and if we can place our problem in one of the scenarios of subsubsection 2.8.1 a great performance leap is to be obtained thanks to the results.

For the 1D $W_\mathcal{I}^p$ space take $\chi_i = L_i$, $i = 0, 1, ..., p$, therefore we are satisfying the requirement (2.54) thanks to (2.98). Then, by means of the properties and definitions in (2.95)-(2.101) we arrive at the following chart to compute the 1D integrals $F_{c;ij}^{rs}$ defined in (2.90).

1. $r = s = 0$ (recall that $F_{c;ij}^{00} = F_{c;ji}^{00}$)
    - $i, j \geq 2$

$$F_{c;ij}^{00} = \int_0^1 L_i(\xi_c)L_j(\xi_c)d\xi_c = \begin{cases} \frac{1}{4(2j-1)^2}\left(\frac{1}{2j+1} + \frac{1}{2j-3}\right) & \text{if } i = j; \\ \frac{-1}{4(2j+3)(2j-1)(2j+1)} & \text{if } i = j+2; \\ 0 & \text{otherwise.} \end{cases} \tag{2.103}$$

- $j = 1$

$$F_{c;i1}^{00} = \begin{cases} 0 & \text{if } i \geq 4; \\ -1/60 & \text{if } i = 3; \\ -1/12 & \text{if } i = 2; \\ 1/3 & \text{if } i = 1. \end{cases} \quad (2.104)$$

- $j = 0$

$$F_{c;i0}^{00} = \begin{cases} (-1)^i F_{c;i1}^{00} & \text{if } i \geq 2; \\ 1/6 & \text{if } i = 1; \\ 1/3 & \text{if } i = 0. \end{cases} \quad (2.105)$$

2. $r = 0;\ s = 1$
    - $i = 1$

$$F_{c;1j}^{01} = \begin{cases} 0 & \text{if } j \geq 3; \\ 1/6 & \text{if } j = 2; \\ 1/2 & \text{if } j = 1; \\ -1/2 & \text{if } j = 0. \end{cases} \quad (2.106)$$

- $i = 2$

$$F_{c;2j}^{01} = \begin{cases} 1/30 & \text{if } j = 3; \\ -1/6 & \text{if } j = 1; \\ 1/6 & \text{if } j = 0; \\ 0 & \text{otherwise}. \end{cases} \quad (2.107)$$

- $i \geq 3;\ j \geq 1$

$$F_{c;ij}^{01} = \frac{1}{2(2i-1)} \left[ \delta_{i,j-1} \frac{1}{2i+1} - \delta_{i-1,j} \frac{1}{2i-3} \right]. \quad (2.108)$$

- $i = 0;\ j \geq 1$

$$F_{c;0j}^{01} = \begin{cases} 0 & \text{if } j \geq 3; \\ -1/6 & \text{if } j = 2; \\ 1/2 & \text{if } j = 1. \end{cases} \quad (2.109)$$

- $i \geq 0;\ j = 0$

$$F_{c;i0}^{01} = -F_{c;i1}^{01} \quad (2.110)$$

- Antisymmetry. From above we deduct:

$$F_{c;ij}^{10} = -F_{c;ji}^{10} \text{ if one of the following holds } \begin{cases} i \geq 2, j \geq 0 & \text{or} \\ i \geq 0, j \geq 2 & \text{or} \\ i = 0, j = 1 & . \end{cases} \quad (2.111)$$

3. $r = 1;\ s = 0$

$$F_{c;ij}^{10} = F_{c;ji}^{01} \tag{2.112}$$

4. $r = s = 1$ (recall that $F_{c;ij}^{11} = F_{c;ji}^{11}$)
    - $i, j \geq 1$

$$F_{c;ij}^{11} = \delta_{ij}\frac{1}{2i-1} \tag{2.113}$$

    - $j = 0$

$$F_{c;i0}^{11} = \begin{cases} 0 & \text{if } i \geq 2; \\ -1 & \text{if } i = 1; \\ 1 & \text{if } i = 0. \end{cases} \tag{2.114}$$

# 3 Applications and results

As stated earlier, the main motivation for proposing these integration algorithms has been to use them within the computational implementation of DPG finite element methods, as they require the computation of the Gram matrix for a high-order discrete test space for each element. Next, we summarize the principles of DPG to illustrate the benefit of applying a fast integration technique in its implementation.

The DPG methodology is a family of finite element techniques that are applicable with several variational formulations [6], a growing number of physical applications [24, 13, 13, 8, 15, 23, 16, 11], and holds desirable mathematical properties like getting "automatic" stability in any norm [5], being able to conformingly handle general polygonal meshes [25] and having a built-in *a posteriori* error estimator [9, 3], among others. In general terms, we can describe any DPG method by referring to an abstract variational formulation like the following one: find $u \in \mathcal{U}$ such that

$$b(u, v) = \ell(v) \quad \forall v \in \mathcal{V}, \tag{3.1}$$

where $\mathcal{U}, \mathcal{V}$ are Hilbert spaces of measurable functions defined over a domain $\Omega$, $b(\cdot, \cdot) : \mathcal{U} \times \mathcal{V} \to \mathbb{R}$ (or $\mathbb{C}$) is a continuous bilinear (or sesquilinear) form, and $\ell : \mathcal{V} \to \mathbb{R}$ (or $\mathbb{C}$) is a continuous linear (or antilinear) functional on $\mathcal{V}$. Suppose that $\mathcal{T}$ is an open partition of our domain $\Omega$, which coincides with the geometry discretization (mesh) used in the finite element method context. Define the broken test space as $\mathcal{V}(\mathcal{T}) := \{f : \Omega \to \mathbb{R}(\text{or } \mathbb{C}, \text{or } \mathbb{R}^3, ...)$ measurable, such that $f|_\mathcal{K} \in \mathcal{V}|_\mathcal{K}, \mathcal{K} \in \mathcal{T}\}$. If testing with functions from the latter space, an additional set of unknowns is necessarily added. Let $\partial\mathcal{T}$ be the set of mesh interfaces, or mesh "skeleton"; let $\mathcal{W}(\partial\mathcal{T})$ be a space of functions defined on $\partial\mathcal{T}$. The resulting problem is, find $u \in \mathcal{U}, \widetilde{u} \in \mathcal{W}(\partial\mathcal{T})$ such that

$$b(u, v) + \widetilde{b}(\widetilde{u}, v) = \ell(v) \quad \forall v \in \mathcal{V}(\mathcal{T}), \tag{3.2}$$

where $\widetilde{b}(\cdot, \cdot)$ is a new bilinear functional acting on elements of $\mathcal{W}(\partial\mathcal{T})$ paired with traces of $\mathcal{V}(\mathcal{T})$ on $\partial\mathcal{T}$. When going to the discrete (finite-dimensional) problem, we'd like to find $u^h \in \mathcal{U}^h, \widetilde{u}^h \in \mathcal{W}^h(\partial\mathcal{T})$, and for the sake of stability, an appropriate finite-dimensional subspace of $\mathcal{V}(\mathcal{T})$ must be picked. DPG's theoretical derivation leads to an exact way of performing this, but it involves the inversion of the Riesz map of $\mathcal{V}(\mathcal{T})$, which in most cases is impossible. Then we limit ourselves to working with an enriched test space, satisfying $\dim \mathcal{V}^r(\mathcal{T}) > \dim \mathcal{U}^h + \dim \mathcal{W}^h(\partial\mathcal{T})$. When implementing this, a convenient way to enrich the test space is to use polynomial spaces of nominal order $p + \Delta p$, with certain increment $\Delta p > 0$ with respect to the nominal polynomial order $p > 0$ of the exact sequence associated to the discrete trial spaces (i.e. $\mathcal{U}^h$ and $\mathcal{W}^h(\partial \mathcal{T})$). Once this enriched test space is formed, the final numerical problem of DPG becomes: find $\psi \in \mathcal{V}^r(\mathcal{T}), u^h \in \mathcal{U}^h, \widetilde{u}^h \in \mathcal{W}^h(\partial \mathcal{T})$ such that

$$\begin{cases} (\psi, v)_{\mathcal{V}^r(\mathcal{T})} - b(u^h, v) - \widetilde{b}(\widetilde{u}^h, v) &= -\ell(v) \quad \forall v \in \mathcal{V}^r(\mathcal{T}) \\ b(\delta u, \psi) &= 0 \quad \forall \delta u \in \mathcal{U}^h \\ \widetilde{b}(\delta \widetilde{u}, \psi) &= 0 \quad \forall \delta \widetilde{u} \in \mathcal{W}^h(\partial \mathcal{T}), \end{cases} \quad (3.3)$$

where $(\cdot, \cdot)_{\mathcal{V}^r(\mathcal{T})}$ is the inner product of the test space. Let $\{u_j\}_{j=1}^n$ be a basis for the first component of the trial space, with $n = \dim \mathcal{U}^h$; $\{\widetilde{u}_j\}_{j=n+1}^{n+\widetilde{n}}$ be a basis for the second part of the trial space, with $\widetilde{n} = \dim \mathcal{W}^h(\partial \mathcal{T})$; and $\{\phi_i\}_{i=1}^m$ be a basis for $\mathcal{V}^r(\mathcal{T})$, where $m$ is its space dimension. With these bases we represent the unknowns as follows:

$$u^h = \sum_{j=1}^n u_j \mathsf{u}_j, \quad (3.4)$$

$$\widetilde{u}^h = \sum_{j=n+1}^{n+\widetilde{n}} \widetilde{u}_j \mathsf{w}_j, \quad (3.5)$$

$$\psi = \sum_{j=1}^m \phi_j \mathsf{s}_j. \quad (3.6)$$

Following this, we can arrive at the augmented linear system shown in the Introduction of this paper. Here, we restate it since the context may now demonstrate that (3.3) is equivalent to this system (in the real case): find $\mathsf{s} \in \mathbb{R}^m, \mathsf{u} \in \mathbb{R}^n, \mathsf{w} \in \mathbb{R}^{\widetilde{n}}$ such that:

$$\begin{cases} \mathsf{G}\mathsf{s} - \mathsf{B}\mathsf{u} - \widetilde{\mathsf{B}}\mathsf{w} &= -\mathsf{l} \\ \mathsf{B}^\mathsf{T}\mathsf{s} &= 0 \\ \widetilde{\mathsf{B}}^\mathsf{T}\mathsf{s} &= 0, \end{cases} \quad (3.7)$$

where $\mathsf{G} \in \mathbb{R}^{m \times m}$ is the Gram matrix of $\mathcal{V}^r(\mathcal{T})$, that is, $\mathsf{G}_{ik} = (\phi_i, \phi_k)_{\mathcal{V}^r(\mathcal{T})}$; $\mathsf{B} \in \mathbb{R}^{m \times n}$ is the enriched stiffness matrix from the field unknowns, given by $\mathsf{B}_{ij} = b(u_j, \phi_i)$; $\widetilde{\mathsf{B}} \in \mathbb{R}^{m \times \widetilde{n}}$ is the enriched stiffness matrix from the skeleton unknowns, with $\widetilde{\mathsf{B}}_{ij} = \widetilde{b}(\widetilde{u}_j, \phi_i)$, and $\mathsf{l} \in \mathbb{R}^m$ is the enriched load vector, obtained through $\mathsf{l}_i = \ell(\phi_i)$.

Because the fundamental purpose of stating this equation system is to find the approximate solution $u^h$ using (3.4), and eventually the solution at the interfaces $\widetilde{u}^h$ through (3.5), it is viable to statically condensate $\mathsf{s}$ from (3.7) by bringing into consideration that the Gram matrix is invertible. Then, the new discrete system yields (also repeating from the Introduction),

$$\begin{pmatrix} \mathsf{B}^\mathsf{T}\mathsf{G}^{-1}\mathsf{B} & \mathsf{B}^\mathsf{T}\mathsf{G}^{-1}\widetilde{\mathsf{B}} \\ \widetilde{\mathsf{B}}^\mathsf{T}\mathsf{G}^{-1}\mathsf{B} & \widetilde{\mathsf{B}}^\mathsf{T}\mathsf{G}^{-1}\widetilde{\mathsf{B}} \end{pmatrix} \begin{pmatrix} \mathsf{u} \\ \mathsf{w} \end{pmatrix} = \begin{pmatrix} \mathsf{B}^\mathsf{T}\mathsf{G}^{-1}\mathsf{l} \\ \widetilde{\mathsf{B}}^\mathsf{T}\mathsf{G}^{-1}\mathsf{l} \end{pmatrix}. \quad (3.8)$$

An alternative way of solving the discrete DPG problem is to apply the so-called Discrete Least Squares framework [17], which has been devised to get, depending on the technique therein implemented, either a more efficient assembly process or a better conditioned linear system than when the full system (3.8) is constructed.

Even though there are not results for hexahedral elements, some stability analysis on triangles and tetrahedra have been made [4, 20]. Such studies showed that the higher the dimension of the enriched test space is with respect to the dimension of the trial spaces, the better is the possibility for the practical DPG method of having stability guaranteed. This implies that the Gram matrix in DPG is the largest array whose computation is required at every element of the partition, therefore a technique to accelerate its construction will favor the overall implementation of this finite element methodology. However, the fact of using a broken test space makes G block diagonal, with each block associated to one and only one element of the mesh. This means that G can be assembled and even "inverted" element-wise, and this is per se a big computational relief. We encourage the reader to refer to another publication on the matter to know more about the specifics of the DPG family of methods [7].

To show the effect in time savings that the proposed tensorization-based integration generates, three different PDEs (along with a variational formulation) were picked and we compared the time needed to calculate the Gram matrix using conventional, tensorization, and simplified tensorization algorithms. Using the first two algorithms, it was compared also the time elapsed for the formation of G, B and l altogether. For all the examples, the physical domain was the master hexahedron, and a single-element mesh was considered. In the following examples, the ordering of the loops within the algorithms that was preferred is the alternative one presented in Algorithm 3.

As a final but important remark, the matrices G, B and l obtained in each of the examples below are equal (up to machine precision) independently of the algorithm used, which verifies the proper implementation of the different tensorization-based algorithms. In order to make such a verification step, the resulting values from the conventional integration algorithm were used as the reference.

## 3.1 Primal formulation for the Poisson problem

Poisson's equation is given by
$$-\operatorname{div}(k\nabla u) = r, \tag{3.9}$$
where $k$ is a diffusivity coefficient, and $r$ is a source term. This PDE, in the context of DPG will lead to the broken primal variational formulation, in the terms of (3.2),

$$\begin{aligned}
\mathcal{U} &= H^1(\Omega), \\
\mathcal{W}(\partial\mathcal{T}) &= H^{-1/2}(\partial\mathcal{T}) \\
\mathcal{V} &= H^1(\mathcal{T}), \\
b(u,v) &= (k\nabla u, \nabla v)_{\mathcal{T}}, \\
\widetilde{b}(\widetilde{u},v) &= \langle \widetilde{u}, v\rangle_{\partial\mathcal{T}}, \\
\ell(v) &= (r,v)_{\mathcal{T}},
\end{aligned} \tag{3.10}$$

where $(\cdot,\cdot)_{\mathcal{T}} = \sum_{\mathcal{K}\in\mathcal{T}}(\cdot,\cdot)_{\mathcal{K}}$, $\langle\cdot,\cdot\rangle_{\partial\mathcal{T}}$ is a duality pairing between $H^{-1/2}(\partial\mathcal{T})$ and $H^{1/2}(\partial\mathcal{T})$. Recalling the Piola maps defined above, the main discrete element subspaces are $\mathcal{U}^h = W^{p_0} = T^{\text{grad}}\mathcal{Q}^{p_0,p_0,p_0}(\mathcal{I}^3)$ and $\mathcal{V}^r = W^{p_r} = T^{\text{grad}}\mathcal{Q}^{p_r,p_r,p_r}(\mathcal{I}^3)$, where $p_r = p_0 + \Delta p$. In order to compute the enriched stiffness matrix, and the enriched load vector a new portion of code must be added to algorithms 4 and 5. The procedure carried out to compute these additional arrays so far is not optimal, so in a later work the performance of this task can be improved. Moreover, the computation of $\widetilde{\mathsf{B}}$ is done in the conventional way at all times, hence we leave it out of the computing time measurements.

Table 1 shows the average wall clock time of 50 runs of computing the $H^1$ Gram matrix for this problem, for different values of $p_r$. Besides the conventional and general tensorization-based integration algorithms, the results with the second simplified tensorization case ($\sim$ algorithm 11) are shown, making use of the formulae provided in Subsection 2.8.2. This plot is made with $p_r + 1$ as the abscissae because this is the parameter with respect to which the computational cost estimates were presented. The expected trend lines accompany the three data plots in Figure 1.

| $p_0$ | $\Delta p$ | $p_r$ | **Conventional** | **Tensorization** | **Simpl. Tensorization** |
|---|---|---|---|---|---|
| 1 | 1 | 2 | 1.76E-03 | 1.64E-03 | 1.67E-03 |
| 1 | 2 | 3 | 3.14E-03 | 1.95E-03 | 1.89E-03 |
| 2 | 2 | 4 | 1.20E-02 | 3.31E-03 | 2.85E-03 |
| 3 | 2 | 5 | 4.31E-02 | 6.70E-03 | 4.96E-03 |
| 4 | 2 | 6 | 0.163 | 1.49E-02 | 9.45E-03 |
| 5 | 2 | 7 | 0.522 | 3.05E-02 | 1.80E-02 |

Table 1: Average computation time (seconds) of $G^{\text{grad}}$ in the Primal Poisson DPG implementation, for different polynomial orders and three integration algorithms

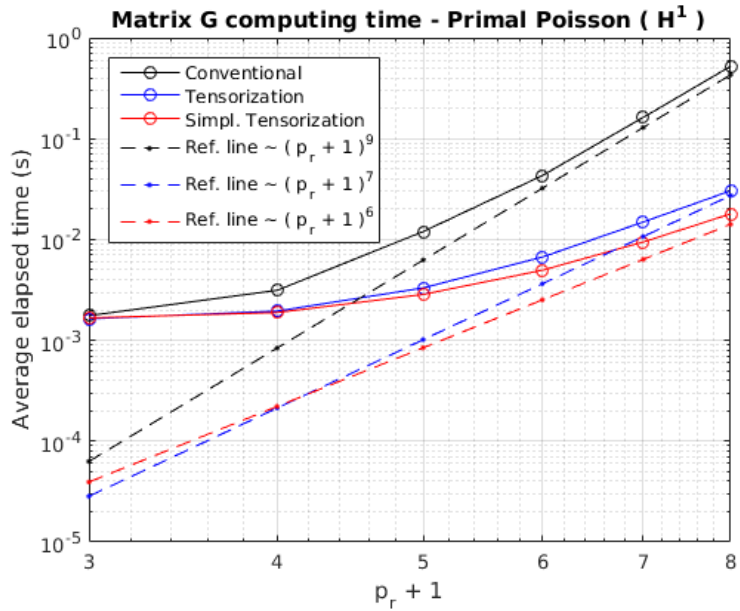

Figure 3.1: Average computation time (seconds) of $G^{\text{grad}}$ in the Primal Poisson DPG implementation, for different polynomial orders and three integration algorithms with respect to $p_r + 1$

Additionally, Table2 and Figure 3.2 show the average time for computing G, B and l within the *elem* subroutine, having fixed $\Delta p = 2$. Since B and l are constructed using a conventional integration approach, this result intends to inform a closer estimate of the actual time savings when using tensorization-based algorithms in a DPG computational solution. Here, it is more useful to see how expensive the process is for a given trial space order, then we plot with respect to $p_0$.

| $p_0$ | $\Delta p$ | $p_r$ | **Conventional** | **Tensorization** |
|---|---|---|---|---|
| 1 | 2 | 3 | 3.55E-03 | 2.41E-03 |
| 2 | 2 | 4 | 1.60E-02 | 7.20E-03 |
| 3 | 2 | 5 | 6.62E-02 | 2.81E-02 |
| 4 | 2 | 6 | 0.294 | 0.113 |
| 5 | 2 | 7 | 1.01 | 0.404 |

Table 2: Average computation time (seconds) of $G^{\mathrm{grad}}$, B, l, in the Primal Poisson DPG implementation, for different polynomial orders and two integration algorithms

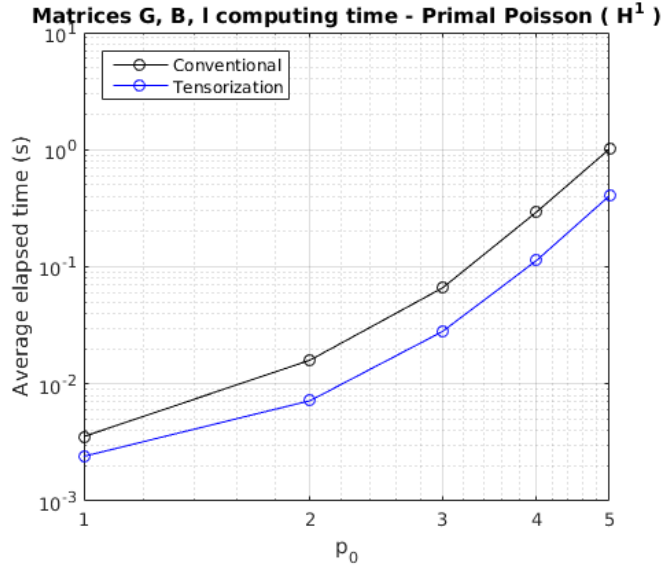

Figure 3.2: Average computation time (seconds) of $G^{\mathrm{grad}}$, B, l, in the Primal Poisson DPG implementation, for different polynomial orders and three integration algorithms with respect to $p_0$, with fixed $\Delta p = 2$

## 3.2 Primal formulation for Maxwell's equations

A second illustrative problem is the time-harmonic Maxwell's wave equation,

$$\operatorname{curl} \mu^{-1} \operatorname{curl} E - \omega^2 \epsilon E = \mathrm{i}\omega J, \tag{3.11}$$

where $E$ is the unknown electric field, $J$ is the imposed electric current, $\omega$ is the frequency, and $\mu, \epsilon$ are electromagnetic properties. In terms of the abstract formulation (3.2), the primal version of the broken variational form for (3.11) is

$$\begin{aligned}
\mathcal{U} &= H(\operatorname{curl}, \Omega), \\
\mathcal{W}(\partial \mathcal{T}) &= H^{-1/2}(\operatorname{div}, \partial \mathcal{T}) \\
\mathcal{V} &= H(\operatorname{curl}, \mathcal{T}), \\
b(E, F) &= (\mu^{-1} \operatorname{curl} E, \operatorname{curl} F)_{\mathcal{T}} - \omega^2 (\epsilon E, F)_{\mathcal{T}}, \\
\widetilde{b}(\widetilde{E}, F) &= -\mathrm{i}\omega \langle \widetilde{E}, F \rangle_{\partial \mathcal{T}}, \\
\ell(F) &= \mathrm{i}\omega (J, F)_{\mathcal{T}},
\end{aligned} \tag{3.12}$$

where $E$ and $\widetilde{E}$ are the trial variables, $F$ is the test variable, and $\langle \cdot, \cdot \rangle_{\partial \mathcal{T}}$ in this case represents the duality pairing between trace spaces $H^{-1/2}(\text{div}, \partial \mathcal{T})$ and $H^{-1/2}(\text{curl}, \partial \mathcal{T})$ [4].

The discrete enriched test subspace is $Q^{p_r} = T^{\text{curl}}(\mathcal{Q}^{p_r-1, p_r, p_r}(\mathcal{I}^3) \times \mathcal{Q}^{p_r, p_r-1, p_r}(\mathcal{I}^3) \times \mathcal{Q}^{p_r, p_r, p_r-1}(\mathcal{I}^3))$. The Gram matrix that corresponds to this problem is the one for $H(\text{curl})$, hereby algorithms 8 and 9 are going to be compared, besides the simplified version of the tensorization-based one.

Similarly to the previous case, we present results for the construction of the Gram matrix in Table 3 and Figure 3.3. In this case, the time averages correspond to 20 runs. When taking into account also the integration of B and l the results are those shown in Table 4 and in Figure 3.4.

| $p_0$ | $\Delta p$ | $p_r$ | Conventional | Tensorization | Simpl. Tensorization |
|---|---|---|---|---|---|
| 1 | 1 | 2 | 3.09E-02 | 3.41E-02 | 3.17E-02 |
| 1 | 2 | 3 | 4.22E-02 | 3.98E-02 | 3.32E-02 |
| 2 | 2 | 4 | 0.129 | 5.39E-02 | 3.67E-02 |
| 3 | 2 | 5 | 0.571 | 9.42E-02 | 5.02E-02 |
| 4 | 2 | 6 | 2.30 | 0.203 | 7.88E-02 |
| 5 | 2 | 7 | 7.82 | 0.464 | 0.132 |

Table 3: Average computation time (seconds) of $G^{\text{curl}}$ in the Primal Maxwell DPG implementation, for different polynomial orders and three integration algorithms

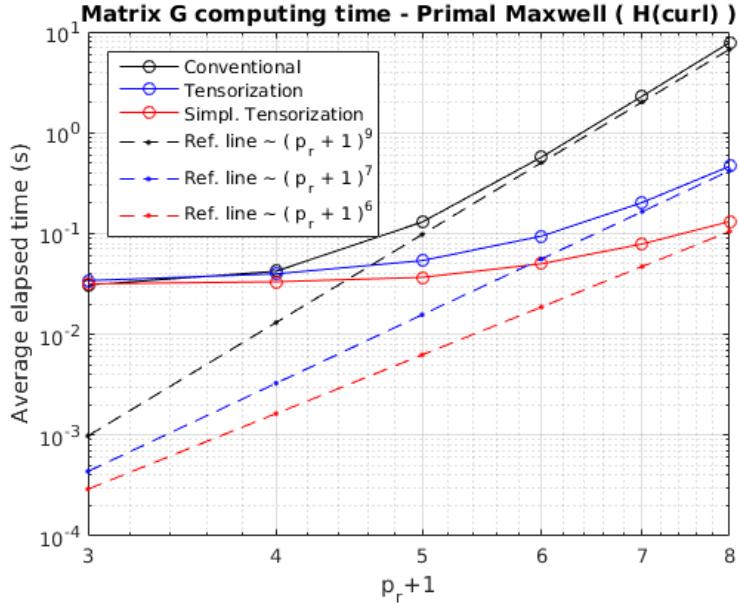

Figure 3.3: Average computation time (seconds) of $G^{\text{curl}}$ in the Primal Maxwell DPG implementation, for different polynomial orders and three integration algorithms with respect to $p_r + 1$

## 3.3 Acoustics and the ultraweak variational formulation

The time-harmonic acoustics system of equations reads as follows:

| $p_0$ | $\Delta p$ | $p_r$ | **Conventional** | **Tensorization** |
|---|---|---|---|---|
| 1 | 2 | 3 | 4.65E-02 | 3.37E-02 |
| 2 | 2 | 4 | 0.167 | 7.75E-02 |
| 3 | 2 | 5 | 0.889 | 0.340 |
| 4 | 2 | 6 | 4.02 | 1.57 |
| 5 | 2 | 7 | 15.36 | 6.02 |

Table 4: Average computation time (seconds) of $G^{\text{curl}}$, B, l in the Primal Maxwell DPG implementation, for different polynomial orders and two integration algorithms

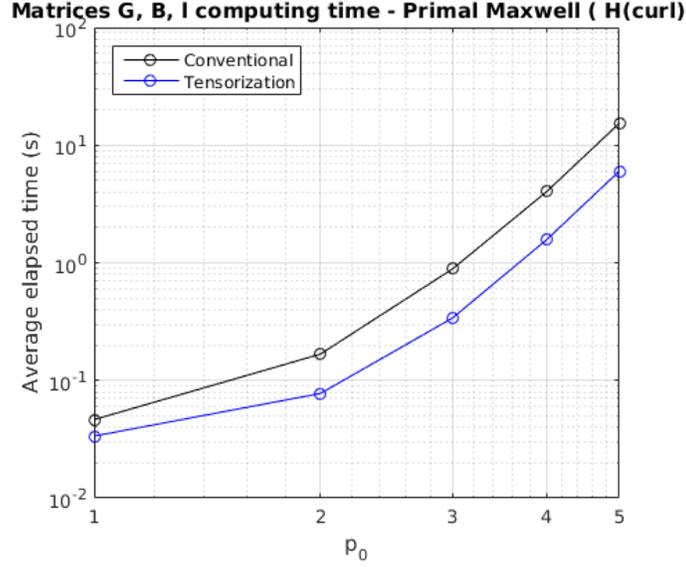

Figure 3.4: Average computation time (seconds) of $G^{\text{curl}}$, B, l, in the Primal Maxwell DPG implementation, for different polynomial orders and three integration algorithms with respect to $p_0$, with fixed $\Delta p = 2$

$$\begin{cases} i\omega u + \nabla p = 0, \\ i\omega p + \text{div}\, p = f, \end{cases} \tag{3.13}$$

where $u$ is the velocity, $p$ the pressure, and $\omega$ the angular frequency. As a first order system, each equation can be tested independently (with broken test functions), and after applying integration by parts we get the so-called ultraweak variational formulation [23]:

$$\begin{aligned}
\mathcal{U} &= L^2(\Omega) \times \boldsymbol{L}^3(\Omega), \\
\mathcal{W}(\partial\mathcal{T}) &= H^{1/2}(\text{div}, \partial\mathcal{T}) \times H^{-1/2}(\text{div}, \partial\mathcal{T}) \\
\mathcal{V} &= H^1(\mathcal{T}) \times H(\text{div}, \mathcal{T}), \\
b((p,u),(q,v)) &= (i\omega u, v)_\mathcal{T} - (p, \text{div}\, v)_\mathcal{T} + (i\omega p, q)_\mathcal{T} - (u, \nabla q)_\mathcal{T}, \\
\widetilde{b}((\widetilde{p}, \widetilde{u}_n),(q,v)) &= \langle \widetilde{p}, v \cdot n \rangle_{\partial\mathcal{T}} + \langle \widetilde{u}_n, q \rangle_{\partial\mathcal{T}}, \\
\ell((q,v)) &= (f,v)_\mathcal{T},
\end{aligned} \tag{3.14}$$

Interestingly, this variational formulation is characterized by having its field trial variables lying in $L^2$

spaces. This immediately implies that there is no requirement on continuity across elements for the functions belonging to those spaces. Other variational formulations, where the trial space needs some kind of continuity, and therefore the handling of orientation and the distinction among the type of shape function (vertex, edge, face, interior) make the tensorization less trivial (see [12]). Unlike those, the ultraweak variational formulation directly allows the correspondence of its discrete trial basis with tensor-product shape functions. As a result, we can apply very similar ideas as above to compute the volume integrals needed for B and l, in addition to G.

For a single element, notice that the stiffness matrix B arising from the bilinear functional in (3.14) can be computed as follows,

$$\mathsf{B}_{IK} = (\mathrm{i}\omega u_{K_R}, v_{I_V})_{\mathcal{K}} - (p_{K_Q}, \operatorname{div} v_{I_V})_{\mathcal{K}} + (\mathrm{i}\omega p_{K_Q}, q_{I_H})_{\mathcal{K}} - (u_{K_R}, \nabla q_{I_H})_{\mathcal{K}}, \qquad (3.15)$$

where $I$ is the product index of $I_H$ and $I_V$, $K$ is the product index of $K_Q$ and $K_R$, $p_{K_Q} \in Y^{p_0}$, $u_{K_R} \in (Y^{p_0})^3$, $q_{I_H} \in W^{p_r}$, and $v_{I_V} \in V^{p_r}$. By mapping to the master hexahedron and using the Piola transformations we obtain:

$$\mathsf{B}_{IK} = \int_{\hat{\mathcal{K}}} \left( \mathrm{i}\omega \hat{u}_{K_R}^\mathsf{T} \mathcal{J} \hat{v}_{I_V} |\mathcal{J}|^{-1} - \hat{p}_{K_Q} \hat{\operatorname{div}} \hat{v}_{I_V} |\mathcal{J}|^{-1} + \mathrm{i}\omega \hat{p}_{K_Q} \hat{q}_{I_H} - \hat{u}_{K_R}^\mathsf{T} \mathcal{J}^{-\mathsf{T}} \hat{\nabla} \hat{q}_{I_H} \right) d^3\boldsymbol{\xi}. \qquad (3.16)$$

For each term of (3.16), we can express the involved shape functions as tensor products of univariate polynomials, and develop similar algorithms as the ones derived above. A similar idea can be applied for the simpler situation of l.

Additionally, the choice of a special norm in the ultraweak formulation is of great benefit regarding the quality of the numerical solution [6]. Such a norm is the so-called *adjoint graph norm*, with a corresponding inner product that can be explicitly calculated as follows,

$$\begin{aligned}((\delta q, \delta v), (q, v))_{\mathcal{V}} =& [\omega^2 + \alpha](\delta q, q)_{\mathcal{T}} + (\nabla \delta q, \nabla q)_{\mathcal{T}} - \mathrm{i}\omega(\operatorname{div} \delta v, q)_{\mathcal{T}} + \mathrm{i}\omega(\delta v, \nabla q)_{\mathcal{T}} \\ &+ \mathrm{i}\omega(\delta q, \operatorname{div} v)_{\mathcal{T}} - \mathrm{i}\omega(\nabla \delta q, v)_{\mathcal{T}} + [\omega^2 + \alpha](\delta v, v)_{\mathcal{T}} + (\operatorname{div} \delta v, \operatorname{div} v)_{\mathcal{T}}, \quad \forall (\delta q, \delta v), (q, v) \in \mathcal{V}.\end{aligned}$$
(3.17)

In (3.17), $\alpha$ is some positive real number that is used as a correction due to the fact of using broken test spaces [23]. Notice that four terms in this inner product were already studied in the $H^1$ and $H(\operatorname{div})$ spaces. The remaining four terms form a Hermitian array, so that we need to implement only two of them. Considering a single element, the third and fourth terms of (3.17) lead to

$$-\mathrm{i}\omega(\operatorname{div} \delta v, q)_{\mathcal{K}} + \mathrm{i}\omega(\delta v, \nabla q)_{\mathcal{K}} = \mathrm{i}\omega \int_{\hat{\mathcal{K}}} \left[ -(\hat{\operatorname{div}}\hat{\delta v})(\hat{q}) + (\hat{\delta v})^\mathsf{T}(\hat{\nabla}\hat{q}) \right] d^3\boldsymbol{\xi}, \qquad (3.18)$$

after implementing the respective Piola maps. Finally, besides using Algorithms 5 and 7 for the four associated terms in (3.17), an analogous tensorization-based algorithm must be built for the two integrals in (3.18).

This problem was solved using both conventional integration and tensorization-based integration, translating the new ideas for the stiffness matrix, the load vector and the Gram matrix produced by the adjoint graph norm, into the code. Table 5 and Figure 3.5 show the results of the average computing time (of 20 runs) of G, B, and l, for several polynomial degrees.

| $p_0$ | $\Delta p$ | $p_r$ | Conventional | Tensorization |
|---|---|---|---|---|
| 1 | 1 | 2 | 1.48E-03 | 7.15E-04 |
| 1 | 2 | 3 | 1.05E-02 | 3.40E-03 |
| 2 | 2 | 4 | 7.46E-02 | 9.62E-03 |
| 3 | 2 | 5 | 0.437 | 2.84E-02 |
| 4 | 2 | 6 | 1.900 | 8.08E-02 |
| 5 | 2 | 7 | 6.985 | 0.214 |
| 6 | 2 | 8 | 21.88 | 0.496 |

Table 5: Average computation time (seconds) of G, B and l in the ultraweak acoustics DPG implementation, for different polynomial orders and two integration algorithms

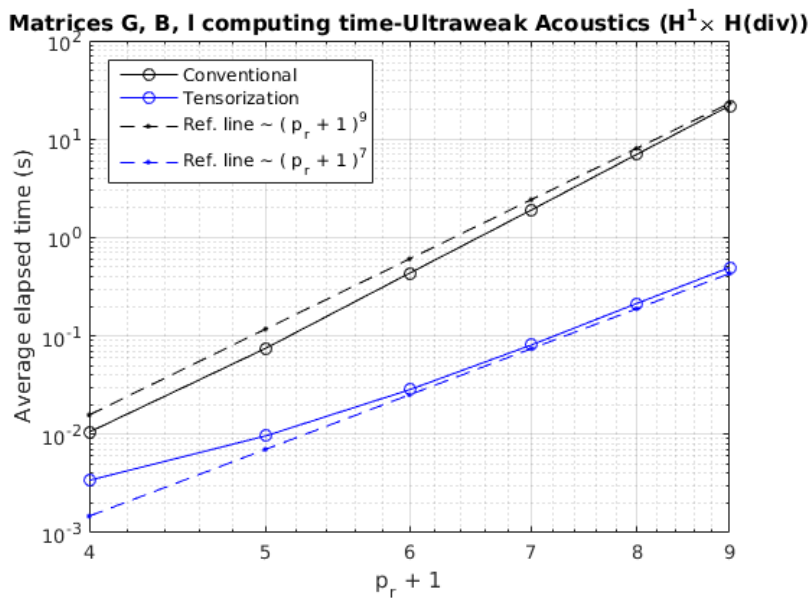

Figure 3.5: Average computation time (seconds) of G, B, l, in the ultraweak acoustics DPG implementation, for different polynomial orders and two integration algorithms, with fixed $\Delta p = 2$

## 3.4 Discussion

The results observed above allow to identify several achievements, while also make clear what aspects need improvement.

Firstly, Figures 3.1 and 3.3 communicate well that the expected asymptotic reference lines approach the numerical experiments results as $p_r$ increases, at least in the conventional and full tensorization-based algorithm. This is more noticeable in the case of Poisson, in which, because of dealing with a $H^1$ Gram matrix and a smaller set auxiliary functions are required (in algorithm 5), then the final accumulation statement quickly begins to dominate the computational complexity as the polynomial order rises. On the other hand, the case of Maxwell's equation required the use of algorithm 9, which has intermediate accumulation statements that are comparable in cost to the final one when $p_r$ is not very high. However, on Figure 3.3 the last three values of polynomial order indicate that the computing time clearly develops at a slower rate than in the conventional case, which indeed looks quite close to the trend line $(p_r + 1)^9$.

The next remarkable result lies in the behavior of the simplified tensorization-based algorithm for computing the Gram matrices of the first two problems. Figures 3.1 and 3.3 show how quickly this modified algorithm delivers a considerably cheaper result, although it approaches its reference line at a much slower pace than the other two. Nevertheless, in both Poisson's and Maxwell's problems, for a value of $p_r$ as low as 4, the time saving for computing the Gram matrix is evident when comparing to both the conventional and the full tensorization-based integration approaches. If we compare our first example to the data reported by Kurtz [18], we achieve quite similar results. In both cases, for $p_r + 1 = 8$ a ratio of about 15 was found (See Figure 3.1). Furthermore, if we contrast that number with the result using the simplified map version, the ratio increases to approximately 30.

We remark, however, that incorporating an algorithm like this makes a FE code less generalized because of the type of element map that is required to be able to use such an approach. Even though it can be recommended to implement this idea when the geometry permits, this implementation has the additional drawback of making the computation of the stiffness matrix and load vector more complicated, therefore it was not taken into account when measuring the computing time of all the matrices.

As a third observation, we see that Figures 3.2 and 3.4 reveal that, when looking at the elapsed time of computing all the matrices, the tensorization-based approach overcomes the performance of the conventional algorithm at all degrees $p_0$. Such a result is also visible in Tables 2 and 4, in which the difference between columns 4 and 5 correspond to the actual time saving when going from the conventional integration scheme to the herein introduced algorithm, for every time the *elem* subroutine is invoked in a DPG code. For the highest values of $p_0$, this saving can be of the order of a second per element in Poisson, or 9 seconds in Maxwell, therefore leading to a big reduction of computation time when a fine mesh is in use.

The advantages of using the tensorization-based approach for numerical integration has, nonetheless, the requirement of finding a procedure of computing B and l that can take more profit of the nesting of integrals developed above for so diverse situations. Clearly, for the ultraweak formulation that was studied, there was a way to extend the ideas of sum-factorization and take profit of them, as noticed in Figure 3.5, where it can be seen that the complexity rates were successfully taken to the desired order. However, for the first two problems, the effect of not thinking through a more efficient algorithm for these matrices can be appreciated in the last points of Figures 3.2 and 3.4, where, despite relatively big savings in the overall computing time of the matrices under study, it is clear that the rate of the tensorization-based algorithm is similar to that of the conventional algorithm. This happens because computing the stiffness matrix B in this way has a computational complexity leading term of $\mathcal{O}((p_0 + 1)^9)$, which will obviously begin dominating the cost of all the integration as $p_0$ grows.

As a closing remark, we can state that having proposed and utilized a tensorization-based algorithm for the integration of non-trivial Gram matrices, like that of $H(\text{curl})$ or the adjoint-graph inner product for $H^1 \times H(\text{div})$, proved to be very useful regarding the enormous time saving (about 8 seconds saved per element when $p_0 = 5$ and $\Delta p = 2$, in Maxwell's equation with the simplified version; or 21 seconds saved in the ultraweak formulation for acoustics with $p_0 = 6$ and $\Delta p = 2$, without simplification), hence this kind of benefit is expected to extend to the rest of energy spaces and inner products frequently used in DPG and other FE methods.

# 4 Conclusions

A complete set of algorithms for fast integration of Gram matrices for $H^1$, $H(\text{curl})$, $H(\text{div})$ and $L^2$ spaces for a general parametric hexahedral element has been presented. The algorithms are based on the sum factorization of tensor-product shape functions. Critical for efficiency is the use of hierarchical shape functions - Legendre polynomials and their integrals - for the exact 1D sequence. The expected reduction of computational complexity from $\mathcal{O}(p^9)$ to $\mathcal{O}(p^7)$ was verified numerically for the $H^1$, $H(\text{div})$ and $H(\text{curl})$ cases. Of special significance was the implementation of a problem with the ultraweak variational formulation, that allowed the full application of the tensorization-based integration to both the stiffness matrix and load vector, unlike the other problems studied, for which only the Gram matrix could be treated under this scheme.

Additionally, when we had the particular case of a simpler element map, the algorithms were modified and the computational complexity was taken down to $\mathcal{O}(p^6)$, reduction that was noticeable in the numerical experiments.

Although in some cases, the time saving in obtaining the Gram matrix is as dramatic as going from 22 seconds to half second (acoustics, $p_r = 8$), and this result should persuade some users to consider the implementation of this set of algorithms, some extra work on this idea seems relevant in order to consolidate the concept of fast integration in DPG thanks to a tensorization-based integration. Such future work may include a specific way of dealing with the stiffness and load arrays when the trial shape functions require compatibility between adjacent elements; extending these ideas to the computation of the interfacial contributions to the stiffness matrix and load vector; the application of these algorithms to more boundary-value problems that have been studied with DPG, so that more test spaces can be tried with this approach, as well as more non conventional norms get to be implemented under this perspective; and finally, propose the algorithms for other types of 2D and 3D elements.

**Acknowledgments.** Jaime Mora has been sponsored by a 2015 Colciencias-Fulbright scholarship, granted by the Government of Colombia and the Fulbright Commission-Colombia. Leszek Demkowicz has been supported by a grant from AFOSR (FA9550-12-1-0484).**Acknowledgments.** Jaime Mora has been sponsored by a 2015 Colciencias-Fulbright scholarship, granted by the Government of Colombia and the Fulbright Commission-Colombia. Leszek Demkowicz has been supported by a grant from AFOSR (FA9550-12-1-0484).